\documentclass{article}

\usepackage{arxiv}
\usepackage[utf8]{inputenc} % allow utf-8 input
\usepackage[T1]{fontenc}    % use 8-bit T1 fonts
\usepackage{hyperref}       % hyperlinks
\usepackage{url}            % simple URL typesetting
\usepackage{booktabs}       % professional-quality tables
\usepackage{amsfonts}       % blackboard math symbols
\usepackage{nicefrac}       % compact symbols for 1/2, etc.
\usepackage{microtype}      % microtypography
\usepackage{lipsum}
\usepackage{graphicx}
\usepackage{caption}
\usepackage{float}
\usepackage{subfigure}
\usepackage{amssymb}
\usepackage{amsmath}
\usepackage{latexsym}
\usepackage{algorithm}
\usepackage{algpseudocode}
\usepackage[mathscr]{euscript}
\usepackage{url}
\usepackage{bm}
\usepackage{textcomp}
\newcommand\norm[1]{\left\lVert#1\right\rVert}
\newcommand*\dif{\mathop{}\!\mathrm{d}}
\usepackage[title]{appendix}
\title{Non-Intrusive Reduced-Order Modeling\ Using Convolutional Autoencoders}

\author{
 Rakesh Halder \\
  Department of Aerospace Engineering\\
  University of Michigan\\
  Ann Arbor, MI 48109 \\
   \And
 Krzysztof J. Fidkowski \\
  Department of Aerospace Engineering\\
  University of Michigan\\
  Ann Arbor, MI 48109 \\
  \And
 Kevin J. Maki \\
  Department of Naval Architecture and Marine Engineering\\
  University of Michigan\\
  Ann Arbor, MI 48109 \\
}
\date{\vspace{-5ex}}

\begin{document}
\maketitle

\begin{abstract}
The use of reduced-order models (ROMs) in physics-based modeling and simulation almost always involves the use of linear reduced basis (RB) methods such as the proper orthogonal decomposition (POD). For some nonlinear problems, linear RB methods perform poorly, failing to provide an efficient subspace for the solution space. The use of nonlinear manifolds for ROMs has gained traction in recent years, showing increased performance for certain nonlinear problems over linear methods. Deep learning has been popular to this end through the use of autoencoders for providing a nonlinear trial manifold for the solution space. In this work, we present a non-intrusive ROM framework for steady-state parameterized partial differential equations (PDEs) that uses convolutional autoencoders (CAEs) to provide a nonlinear solution manifold and is augmented by Gaussian process regression (GPR) to approximate the expansion coefficients of the reduced model. When applied to a numerical example involving the steady incompressible Navier-Stokes equations solving a lid-driven cavity problem, it is shown that the proposed ROM offers greater performance in prediction of full-order states when compared to a popular method employing POD and GPR over a number of ROM dimensions.
\end{abstract}

\keywords{model reduction \and autoencoders \and deep learning \and non-intrusive \and machine learning \and proper orthogonal decomposition}

\section{Introduction}
Physics-based modeling and simulation has become an essential tool in many engineering and science applications, allowing for highly accurate representations of physical systems which may be otherwise difficult to evaluate. Physics-based models consist of a set of governing equations, which are often found in the form of parametrized partial differential equations (PDEs) discretized over a computational domain, where a set of design parameters $\bm{\mu}$ controls properties such as the boundary conditions, geometry of the computational domain, or physical properties. In industrial processes such as design optimization, a large number of designs needs to be evaluated and the accuracy, or fidelity of these models must be high. High-fidelity simulations of large scale models are very computationally intensive, requiring large amounts of memory and computational time. This large computational cost can render many-query processes such as design optimization infeasible. 

The use of reduced-order models (ROMs) is a common approach for drastically lowering this computational cost. ROMs create a surrogate model using training data from a set of computed simulations of the high-fidelity full-order model that allows for accurate, rapid, real-time evaluation of simulations at unseen design parameters. ROMs seek to reduce the dimensionality of full-order models, which contain a large number of degrees of freedom. This involves a compression phase in which a \emph{reduced basis} (RB) of the solution space is obtained, from which accurate approximations of full-order solutions can be obtained through a set of \emph{expansion coefficients}. ROMs consist of two stages: a computationally intensive offline stage where high-fidelity solutions are evaluated to obtain data snapshots and a low-dimensional surrogate model is trained, and an online stage where the surrogate model can be rapidly evaluated to approximate solutions at unseen design parameters. 

A commonly used method for obtaining the reduced basis is the proper orthogonal decomposition (POD) \cite{ BerkoozStructure, Holmes1997LowdimensionalMO}, which utilizes the singular value decomposition (SVD) to obtain a low-rank trial subspace composed of a number of linearly independent basis vectors, a linear combination of which is used to approximate unrealized solutions. Projection-based RB methods~\cite{carlberg2011efficient, CARLBERG2017693} project the physics of the governing equations onto the low-rank trial subspace and solve a low-dimensional version of the full-order model. Although projection-based RB methods have been shown to offer robust performance, their computational cost remains large for certain nonlinear problems that have a non-affine dependence on the inputs for quantities such as residuals in computational fluid dynamics (CFD) models~\cite{carlberg2015adaptive}.  Methods such as the discrete empirical interpolation method (DEIM)~\cite{chaturantabut2010nonlinear} allow for an affine representation of operators to avoid full-order evaluations, although realizing this is often difficult and intractable for some nonlinear problems. 

In \emph{non-intrusive} ROMs, an alternative to projection-based ROMs, the governing physics are only used to generate solution snapshots of the high-fidelity model in the offline stage and are not projected onto a lower dimension in the online stage. A regression model is required to interpolate over the expansion coefficients of the training data to approximate them for unseen design parameters. Some popular interpolation methods include Gaussian process regression (GPR)~\cite{guo2018reduced, DupuisGPR} and neural networks ~\cite{HESTHAVEN201855,JacquierNN}. Although RB methods using POD are widely used and generally offer good performance, they often produce inaccurate results for certain nonlinear problems, such as those dominated by advection, and require ROMs of large dimension to produce results with acceptable accuracy~\cite{algoritmy}. For nonlinear problems characterized by different physical regimes, localized POD subspaces~\cite{amsallem2012nonlinear} are often employed to mitigate this issue.

Whereas POD produces a trial subspace that is linear, recent methods have attempted to compute low-dimensional nonlinear \emph{trial manifolds} that are more adept at handling nonlinear problems. Many recent advances have utilized machine learning and artificial intelligence (AI) methods to this end, which have been at the forefront of massive recent breakthroughs in numerous fields such as computer vision, natural language processing, and recommender systems~\cite{krizhevsky2012imagenet, hirschberg2015advances, ZhangRecommender}. The use of machine learning methods has become ubiquitous in many domains, significantly improving and even beating the performance of existing methods. 

There are many machine learning methods for producing low-dimensional representations of high-dimensional data, many of which do so non-linearly as opposed to POD. Deep learning~\cite{lecun2015deeplearning} approaches have been utilized to develop ROMs that provide efficient nonlinear trial manifolds of physical systems. Convolutional autoencoders (CAEs), a type of neural network, have been used in ROMs and have been shown to outperform POD-based methods~\cite{Lee2020ModelRO,Fresca2021ACD}. Convolutional autoencoders are adept at learning data that are spatially distributed, including the solutions to PDEs discretized over a computational domain. Autoencoder neural networks consist of two parts: an encoder, which maps high-dimensional inputs to a low-dimensional code, and a decoder, which maps the low-dimensional code to an approximation of the high-dimensional input. In the context of ROMs, the code is analogous to the expansion coefficients that map back to the full-order solution space. When using autoencoders for ROMs, the entire network is trained in the offline stage, while only the decoder is used in the online stage for rapid evaluation of unseen solutions. To the best of our knowledge, there have been two attempts at using convolutional autoencoders for non-intrusive ROMs; one utilizing them for vehicle aerodynamic simulation~\cite{mrosek2021variational} and another for natural convection in porous media~\cite{kadeethum2022non}. The first found that using autoencoders only offers a very slight improvement over POD-based methods. The second involves unsteady problems with a limited number of design variables and does not elaborate on the projection errors provided by autoencoders and their relation to interpolation accuracy. In this work, we propose a non-intrusive ROM framework consisting of a nonlinear trial manifold produced by a CAE that is augmented by Gaussian process regression to handle the interpolation of expansion coefficients in the design parameter space. This ROM framework is referred to as CAE-GPR and its performance is compared to that of POD-GPR when applied to a problem which solves the incompressible Navier-Stokes equations over a number of ROM dimensions.

\section{Full-order model}
The full-order model (FOM) in this work is considered to be the solution $\bm{x}(\bm{\mu})$ \unboldmath $\in \mathbb{R}^{N}$ of a state variable in a system governed by a set of steady-state parameterized partial differential equations (PDEs) discretized over a computational domain $\Omega \in \mathbb{R}^{d}$. We consider design parameters \boldmath $\mu$ \unboldmath $\in \bm{\mathcal{D}}$ that define both the computational domain and parameters of the governing equations. Here $\bm{\mathcal{D}} \subseteq \mathbb{R}^{p}$
denotes the parameter space such that \boldmath $x$: \unboldmath $ \bm{\mathcal{D}} \rightarrow \mathbb{R}^{N}$. The set of PDEs governing the FOM is solved numerically over $\Omega$ to generate a solution $\bm{x}(\bm{\mu})$. The computational cost of numerically solving the system increases with its dimension $N$, which is in proportion with the fineness of $\Omega$. Accurate or useful solutions of systems often require large values of $N$, resulting in large computational costs for a single solution. In processes such as design optimization, the need to evaluate the solutions for many different designs in real-time becomes infeasible if numerous FOMs have to be solved. This large computational cost motivates the use of reduced-order models, where a small number of FOMs are solved and used to create a computationally inexpensive surrogate model that can deliver accurate approximations in real time.

\section{Linear reduced basis method}
This section gives an overview of the proper orthogonal decomposition (POD), a popular method for constructing a linear reduced basis, which allows for the construction of full-order solutions as a linear combination of independent basis vectors. The basis vectors are formed from a collection of training solutions over the parameter space. A snapshot matrix is assembled from these training solutions, from which the underlying structure of the solution space can be extracted \cite{BerkoozStructure}. Reduced-order models essentially obtain low-dimensional representations of any solution lying on the solution manifold for which there exists a mapping back to the solution space. In classification, a popular choice for linear dimensionality reduction is principal component analysis (PCA), to which the proper orthogonal decomposition is closely related.
\subsection{Proper orthogonal decomposition}
Linear RB reduced-order models rely on training data obtained from a set of $n$ solution snapshots calculated at chosen design points in the parameter space. A snapshot matrix, $\bm{S} \in \mathbb{R}^{N\times n}$, is assembled
\begin{equation}
\bm{S} \in \mathbb{R}^{N\times n}= [{\bm{x}}^1, {\bm{x}}^2, \cdots , {\bm{x}}^n] = [{\bm{x}}(\bm{\mu}^{1}), {\bm{x}}(\bm{\mu}^{2}), \cdots , {\bm{x}}(\bm{\mu}^{n})].
\end{equation}
Denote by $\bm{\mathcal{M}}$ a subspace of the column space of $\bm{S}$. We assume that $\bm{\mathcal{M}}$ provides a good approximation of the solution manifold for $ \bm{\mu} \in \bm{\bm{\mathcal{D}}}$ if there are a sufficient number of solution snapshots in $\bm{S}$ which correspond to a judiciously chosen subset of design parameters in $\bm{\mathcal{D}}$. $\bm{\mathcal{M}}$ is the span of $k$ orthonormal basis vectors,  $[{\bm{\psi}}^{1}, {\bm{\psi}}^{2}, \cdots, {\bm{\psi}}^{k}] \in \mathbb{R}^{N}$, where $k \ll N$. The basis is chosen such that each solution snapshot $\bm{x}^{i}$ in $\bm{S}$ can be well-approximated as a linear combination of the basis vectors
\begin{equation}
\bm{x}^{i} \approx a^{i}_{1}\bm{\psi}^{1} + a^{i}_{2}\bm{\psi}^{2} \cdots + a^{i}_{k}\bm{\psi}^{k}.
\end{equation}
Where $\bm{a}^{i}$ is the set of basis coefficients, or expansion coefficients, for a given solution snapshot.
The truncated singular value decomposition of \bm{$S$}, contains two orthonormal matrices $\bm{U} \in \mathbb{R}^{N\times n}$ and $\bm{V} \in \mathbb{R}^{n\times n}$, as well as a diagonal matrix $\bm{\Sigma} \in \mathbb{R}^{n\times n}$
\begin{equation}
\bm{S = U\Sigma V^{T}}.
\end{equation}
Here, $\bm{U}$ contains a set of $n$ left singular vectors that form an orthonormal basis for the column space of $\bm{S}$, $\bm{V}$ contains a set of $n$ right singular vectors that form an orthonormal basis for the row space of $\bm{S}$, and diag($\bm{\Sigma}$) $\in \mathbb{R}^{n}= [{\sigma}_1, {\sigma}_2, \cdots, {\sigma}_n]$ contains the singular values corresponding to the singular vectors in descending order, $\sigma_{1} \geq \cdots \geq \sigma_{n} \geq 0$. The first $k$ left singular vectors of $\bm{U}$ are chosen to be the basis vectors forming the POD basis, $\bm{\Psi} \in \mathbb{R}^{N\times k} = [{\bm{\psi}}^{1}, {\bm{\psi}}^{2}, \cdots, {\bm{\psi}}^{k}]$. Often, the singular values associated with the basis vectors decay very quickly and only the first $k$ singular vectors are chosen to form the POD basis to preserve only the most dominant basis vectors. To determine the value of $k$, the relative information content of the subspace is evaluated
\begin{equation}
    E(k) = \dfrac{\sum_{j=1}^{k}\sigma_{j}}{\sum_{j=1}^{n}\sigma_{j}}, 
\end{equation}
and $k$ is chosen such that $E(k) \geq \epsilon$, where $\epsilon \in $ [0,1) is chosen somewhat arbitrarily, usually to a value $\epsilon \geq 0.95$ \cite{MrosekROM}. Using the POD basis, full-order solutions at unseen design parameters ${\bm{x}}(\bm{\mu}^{*})$ can be approximated
\begin{equation}
\bm{x}(\bm{\mu}^{*}) \approx \bm{\Psi}\bm{a}^{*} = a^{*}_{1}\bm{\psi}^{1} + a^{*}_{2}\bm{\psi}^{2} \cdots + a^{*}_{k}\bm{\psi}^{k},
\end{equation}
where $\bm{a}^{*}$ can be estimated through a computational model that takes $\bm{\mu}^{*}$ as an input. 
\subsection{Projection error}
A measure of quality of the POD basis is its ability to reconstruct solution snapshots $\bm{x}^{i}$ in $\bm{S}$ with a high degree of accuracy. We first calculate the projection of $\bm{x}^{i}$ onto $\bm{\Psi}$
\begin{equation}
    \hat{\bm{x}}^{i} = \bm{\Psi}\bm{\Psi}^{T}\bm{x}^{i}.
\end{equation}
A measure of the relative error over all of the solution snapshots in $\bm{S}$ is measured through the quantity
\begin{equation}
    \epsilon_{\text{POD}} = \sum_{i=1}^{n} \dfrac{\norm{\bm{x}^{i} - \hat{\bm{x}}^{i}}^{2}} {\norm{\bm{x}^{i}}^2}.
\end{equation}
The Schmidt-Eckart-Young theorem \cite{Eckart1936TheAO} states that the POD basis consisting of the first $k$ left singular vectors found from the SVD of $\bm{S}$ minimizes this error amongst all orthonormal bases of rank $k$.

\begin{figure}[!t]
\centering
\includegraphics[width=0.7\textwidth]{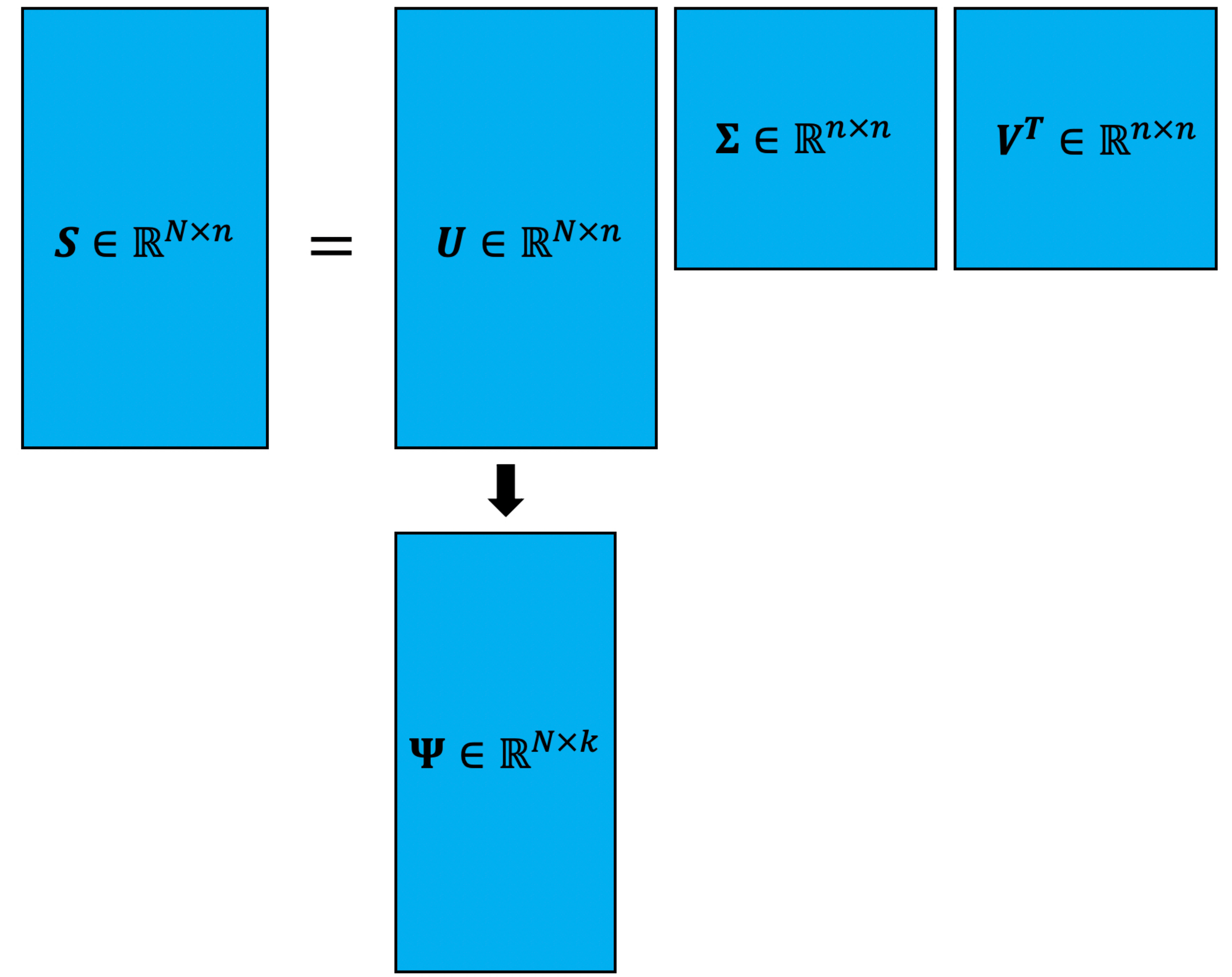}
\caption{Schematic of the proper orthogonal decomposition (POD), where the snapshot matrix $\bm{S}$ is decomposed using the singular value decomposition (SVD) and the POD basis $\bm{\Psi}$ is obtained from $\bm{U}$.}
\end{figure}

\section{Nonlinear manifold construction using convolutional autoencoders}
This section describes the use of deep convolutional autoencoders in constructing a nonlinear trial manifold of simulation data. Unlike a linear reduced basis, which constructs solutions as a linear combination of known and calculated basis vectors, nonlinear manifolds use a mapping function $g(\bm{a})$, which may not be known explicitly, to approximate a mapping between the expansion coefficients and the full-order solution
\begin{equation}
    g(\bm{a}): \mathbb{R}^{k} \rightarrow \mathbb{R}^{N}.
\end{equation}
There are many popular methods for nonlinear dimensionality reduction such as Isomap \cite{tenenbaum2000global} that provide a low-dimensional embedding of high-dimensional data. Nonlinear methods have been shown to offer better performance in classification tasks when compared to linear methods such as PCA \cite{lee2008investigating}, highlighting the advantage of using nonlinear methods to create low-dimensional representations of data. However, most nonlinear dimensionality reduction methods do not provide a mapping back to the high-dimensional solution space which restricts their use in ROMs. Deep convolutional autoencoders, which do provide an approximate mapping $g(\bm{a})$, have been utilized in projection-based ROMs where they have been shown to outperform POD-based methods \cite{Lee2020ModelRO}. Non-intrusive methods have also used deep CAEs to construct nonlinear manifolds that efficiently learn the dynamics of physical systems \cite{gonzalez2018deep}. 

\subsection{Artifical Neural Networks}
\begin{figure}[!t]
\centering
\includegraphics[scale=.375]{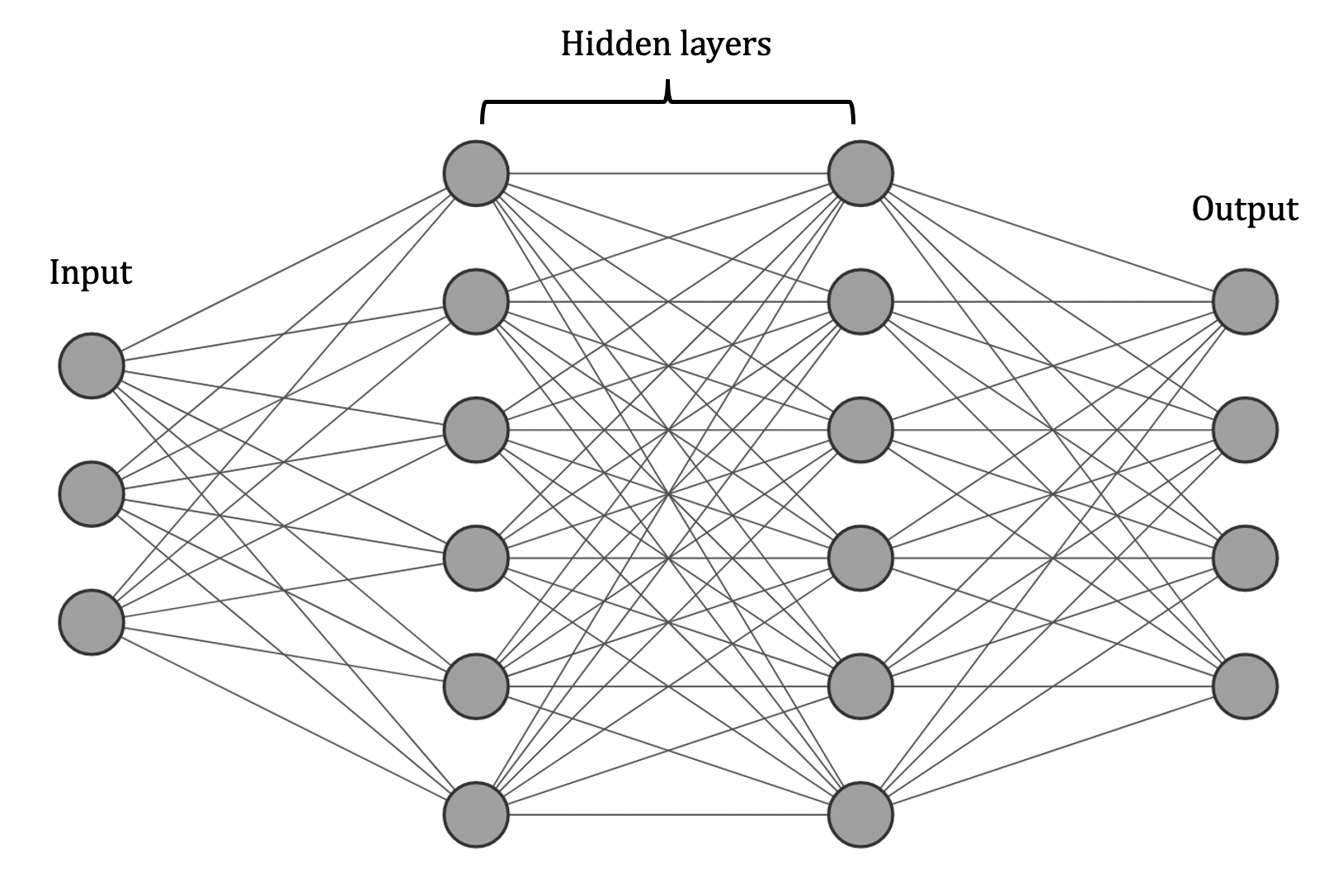}
\caption{Architecture of a multilayer perceptron with a 3-dimensional input, six neurons in two fully connected layers, and four neurons in the output layer.}
\label{fig:ann_architecture}
\end{figure}

An artificial neural network (ANN) is a computational model that is able to learn from a training data set $\bm{\mathcal{T}} = \{\bm{X},\bm{Y}\} $, where $\bm{X}$ and $\bm{Y}$ refer to the inputs and outputs respectively. ANNs are inspired by biological neural networks existing in animal brains \cite{Kriesel2007NeuralNetworks}. ANNs are widely used and versatile models for regression and classification problems. Feedforward neural networks are a type of ANN in which information always propagates in only one direction, creating a direct mapping between inputs and outputs. Feedforward neural networks are composed of an input layer, a number of hidden layers, and an output layer. Neurons comprising these layers are associated with weights and biases, trainable parameters which are optimized during the model training stages.
Figure~\ref{fig:ann_architecture} shows the architecture of a simple feedforward neural network with an input layer, two hidden layers, and an output layer. There are connections between each possible pair of neurons between layers, with each connection carrying a weight term and each neuron carrying a bias term with the exception of those in the input layer. Such a network is referred to as fully connected, or a multilayer perceptron (MLP). Hidden layers in MLPs are also referred to as fully connected layers. Each hidden layer state $\bm{h}_{j}$ is computed from the state in the previous layer,  $\bm{h}_{j-1}$, along with its weights $\bm{W}_{j}$ and biases $\bm{b}_{j}$ as well as an activation function $\phi(x)$

\begin{equation}
    \bm{h}_{j} = \phi\left(\bm{W}_{j}\bm{h}_{j-1}+\bm{b}_{j}\right).
\end{equation}
The role of activation functions is to introduce nonlinearities into the model, allowing for complex functional relationships to arise. In addition to being able to learn the training data well, neural networks should provide reasonable accuracy for unknown data of the same class, a property referred to as generalization \cite{Kriesel2007NeuralNetworks}. A commonly used activation function is the rectified linear unit (ReLU) \cite{nair2010rectified}, which has been shown to offer better performance and ability to generalize when compared to other common activation functions \cite{DahlRELU,Zeiler2013OnRL}

\begin{equation}
    \phi(x) = \text{max}(0,x) = 
    \begin{cases}
    x,& \text{if } x\geq 0\\
    0,              & \text{if } x < 0
\end{cases}.
\end{equation}

For inputs less than 0, the ReLU activation function returns a valuation and gradient of zero, effectively rendering certain neurons inactive. This can be problematic for network training if a large percentage of neurons exhibit this behavior and is commonly referred to as the dying ReLU problem. The leaky ReLU~\cite{xu2015empirical} activation function mitigates this issue, by incorporating a small positive constant $\alpha$ for negative inputs.
\begin{equation}
    \phi(x) =
    \begin{cases}
    x,& \text{if } x\geq 0\\
    \alpha x,              & \text{if } x < 0
\end{cases}.
\end{equation}

\subsection{Training neural networks}
\subsubsection{Backpropagation}
Feedforward networks are trained using a differentiable loss function, $\mathcal{L}\left(\bm{\mathcal{T}}, \left( \bm{W},\bm{b}\right) \right)$, which calculates a measure of error between the state found in the output layer and the correct output values from the training data. The loss function serves as an objective function in an optimization problem, where its gradients with respect to the weights and biases are calculated through backpropagation \cite{Rumelhart:1986we}, an algorithm utilizing automatic differentiation. Common optimizers used in training neural networks include stochastic gradient descent (SGD) and Adam \cite{Adam}. Optimizers perform a number of training epochs over $\bm{\mathcal{T}}$ in an attempt to minimize the loss function. The weights and biases update at the end of epoch $n$ according to
\begin{equation}
    \left(\bm{W}^{n+1},\bm{b}^{n+1}\right) =  \left(\bm{W}^{n},\bm{b}^{n}\right) - \eta \mathcal{G}\left(\frac{\partial \mathcal{L}\left(\bm{\mathcal{T}},\left(\bm{W}^{n},\bm{b}^{n}\right) \right)}{\partial \left(\bm{W}^{n},\bm{b}^{n}\right)}\right),
\end{equation}
where $\eta$ is the learning rate, a hyperparameter controlling the optimizer's step size and $\mathcal{G}$ is a function of the loss function's gradient dependent upon the chosen optimizer. The gradient of the loss function can be calculated using a single training sample as it is when using SGD, using the average gradient of the entire training set, or by using averages of a number of randomly selected mini-batches from $\bm{\mathcal{T}}$. Using mini-batches when training neural networks has been shown to improve the ability to generalize in addition to providing stable convergence \cite{masters2018revisiting}. The mini-batch size $b$ is chosen based on the size of $\bm{\mathcal{T}}$ to strike a balance between performance and computational cost.

The predictive performance of a neural network initially increases with the number of training epochs but starts to stall and then decrease as the network parameters become overly tuned towards the training data and fail to generalize, a problem referred to as overfitting . Regularization methods \cite{YingOverfitting} exist to prevent overfitting. One method is to use early stopping, where the loss on a validation data set $\bm{\mathcal{V}}$ is monitored during training. If $\mathcal{L}\left(\bm{\mathcal{V}},\left(\bm{w},\bm{b}\right)\right)$ fails to drop for a prescribed number of epochs, training is stopped.

The initial set of weights and biases that are used can also effect the final performance of a neural network. A commonly used weight initialization scheme for layers using ReLU activation functions is the He normal~\cite{he2015delving} initializer, which samples weights from a normal distribution centered around 0. It is a common practice to initialize the biases in each layer to 0. There is no standard and accepted approach to choosing the number of hidden layers and the number of nodes in each layer when designing multilayer perceptrons. An optimal choice depends upon a number of factors, including the number of training samples, the dimensionality of the inputs and outputs, the choice of activation functions, and the complexity of the function which is being approximated. The number of hidden layers and nodes to use is often found through a trial and error approach involving model validation techniques such as cross-validation. In general, the total number of trainable parameters in a network is directly related to its capacity to learn functions. Neural networks become deeper as more hidden layers are added. However, network configurations with a large number of trainable parameters tend to overfit to the training data and fail to generalize unless regularization techniques are used. In addition, large networks are computationally expensive to train. In spite of these downfalls, deeper network architectures have become increasingly popular for complex learning tasks in multiple domains as they offer better performance~\cite{lecun2015deeplearning}. Although more than two hidden layers are not required for many learning tasks, some functions are not adequately approximated by networks containing two hidden layers and using deeper networks can drastically improve performance~\cite{telgarsky2016benefits, pmlr-v49-eldan16}.

\subsubsection{Data Normalization}
Similar to many other machine learning algorithms, neural networks often require that the training data be normalized in order to ensure adequate performance \cite{Sola1997ImportanceOI}. Data normalization allows the optimizer to learn the optimal network parameters at a much faster rate. One way to normalize the training data is to apply min-max scaling to each feature in the data matrix $\bm{D}$ containing either the inputs or outputs
\begin{equation}
d' = \frac{d - \text{min}(\bm{d}_{j})}{\text{max}(\bm{d}_{j})-\text{min}(\bm{d}_{j})},
\end{equation}
where $j$ is the feature index. Min-max scaling results in the data being transformed into the range [0,1]. After training, new input data are also normalized while an inverse transformation is applied to predicted outputs. 

\subsection{Autoencoders}
\begin{figure}[!t]
\centering
\includegraphics[width=0.6\textwidth]{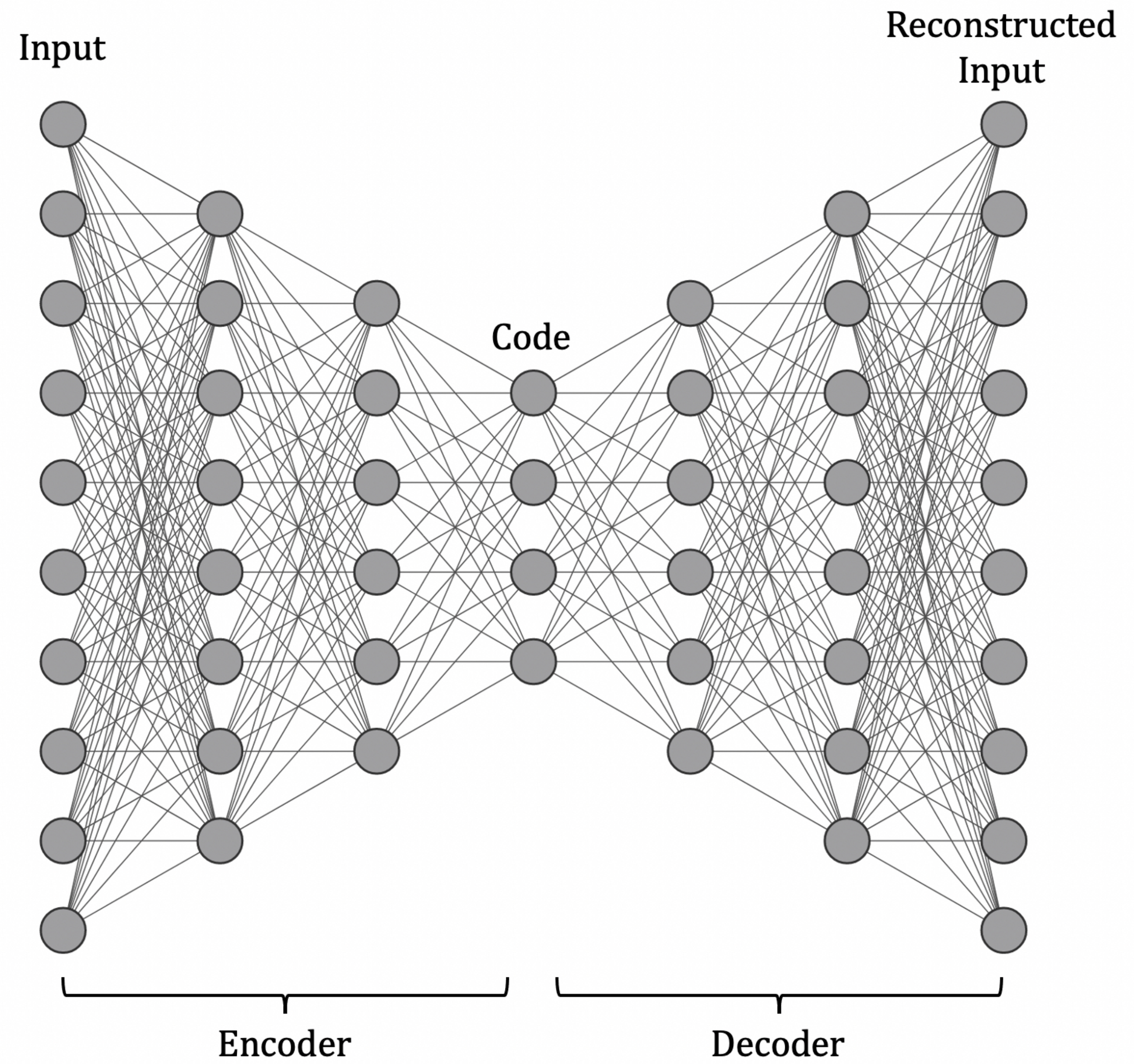}
\caption{Architecture of a symmetric MLP autoencoder with two fully connected layers between the input/output and code.}
\label{fig:autoencoder_architecture}
\end{figure}

Autoencoders are a type of feedforward neural network that aim to learn to reconstruct inputs in the output layer, $g: \bm{x} \rightarrow \bm{\hat{x}}$ where $\bm{x} \approx \bm{\hat{x}}$. Autoencoders use an architecture composed of two individual feedforward neural networks. The encoder $g_{\text{enc}}: \mathbb{R}^{N} \rightarrow \mathbb{R}^{k}$ where $k \ll N$ maps a high-dimensional input $\bm{x}$ into the low-dimensional code $\bm{a}$. The decoder $g_{\text{dec}}: \mathbb{R}^{k} \rightarrow \mathbb{R}^{N}$ maps the code back to an approximation of the high-dimensional input  $\bm{\hat{x}}$. The combination of the two results in
\begin{equation}
    g: \bm{\hat{x}} = g_{\text{dec}} \circ g_{\text{enc}}(\bm{x}). 
\end{equation}
Autoencoders have been shown to provide robust low-dimensional representations of high dimensional data \cite{HintonSalakhutdinov2006b}. Once an autoencoder is sufficiently trained and $g(\bm{x}) \approx \bm{x}$ for all inputs over $\bm{\mathcal{T}}$, the corresponding low-dimensional codes can be passed to the decoder $g_{\text{dec}}(\bm{a})$ to obtain accurate approximations $\bm{\hat{x}}$ for all data in $\bm{\mathcal{T}}$. States existing outside of the training set $\bm{x}^{*}$ can also be well-approximated if a good approximation of the low-dimensional code $\bm{a}^{*}$ can be found. In the context of ROMs, the code is equivalent to the set of expansion coefficients that map from a low-dimensional representation to the high-dimensional full-order solution. Similarly, the projection of a full-order solution onto the nonlinear manifold provided by the autoencoder is given by $\bm{\hat{x}}.$ Training is conducted on the combination of the encoder and decoder, while after training the encoder is often no longer useful and only the decoder is used. Figure~\ref{fig:autoencoder_architecture} shows a sample architecture of a symmetric MLP autoencoder with two hidden layers between the input/output layers and code. Since MLP autoencoders are fully connected, the total number of trainable parameters in the network can grow very large when the dimension of the input, $N$, is high. As the number of trainable parameters increases, the amount of training data required to sufficiently train the network to make reasonably accurate predictions also grows large. This is contrary to the objective of model reduction, which aims to make predictions using a limited amount of training data. 

\subsubsection{Convolutional autoencoders}
\begin{figure}[!t]
\centering
\includegraphics[width=0.95\textwidth]{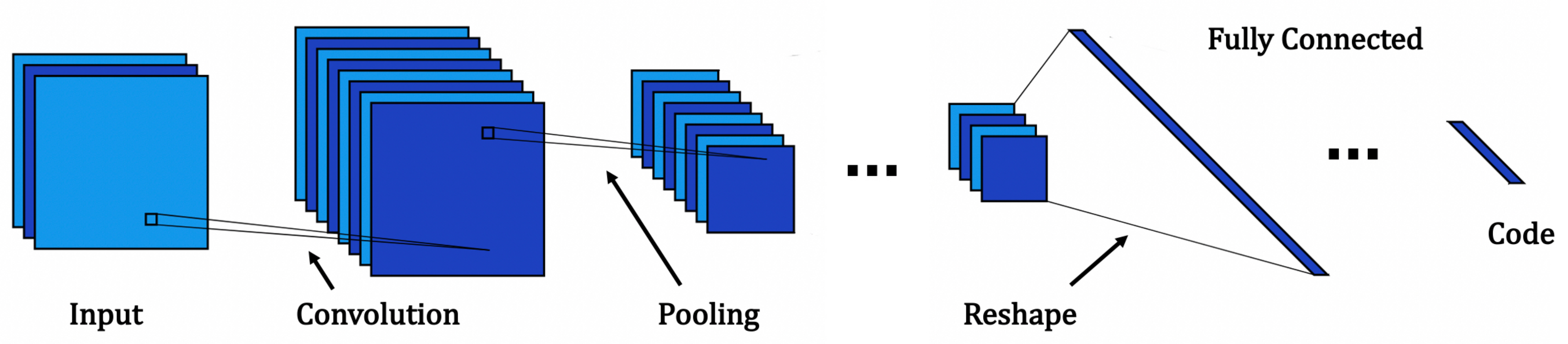}
\caption{Architecture of the encoder of a convolutional autoencoder (CAE) consisting of convolutional, pooling, and fully connected layers.}
\label{fig:cae_encoder}
\end{figure}

There exist neural network architectures that make use of \emph{parameter sharing}, where rather than weight combinations existing for each pair of neurons between layers, multiple neurons share a single weight. \emph{Convolutional autoencoders} effectively implement parameter sharing to limit the total number of trainable parameters in the network. This is done through the use of convolutional layers, which provide feature maps of input data that are spatially arranged \cite{lecun1999object}. Convolutional layers use a number of filters to convolve over spatially distributed input data, with each filter having its own set of weights. Pooling layers are also used in convolutional networks to summarize the features in input data through operations including averaging and maximization. Convolutional layers are widely used in the field of computer vision, dealing with spatially distributed data such as images \cite{simonyan2014very, krizhevsky2012imagenet}. CAEs can also be a useful tool for states that arise from numerically solving discretized PDEs as they tend to be spatially distributed. Data with multiple states, i.e. components of velocity or levels of red, green, and blue in images, can also be handled well by CAEs through the use of a number of input channels. More details on convolutional layers can be found in a work by Dumoulin and Visin~\cite{dumoulin2016guide}. A combination of convolutional, pooling, and fully connected layers is used to construct CAEs, as shown in a schematic of an encoder section of a CAE in Figure~\ref{fig:cae_encoder}. Spatially distributed data arising from the solutions of discretized PDEs often vary smoothly through the computational domain. CAEs are highly adept at handling data that are naturally spatially distributed by learning spatially invariant features, allowing them to outperform other neural network architectures.  \cite{lecun_doc,Goodfellow-et-al-2016}. 
The input and output layers of CAEs usually consist of 2-dimensional (2D) states in each channel. Training data must be reshaped before being input into the network through the use of a reshape operator 
\begin{equation}
    \bm{R}: \mathbb{R}^{N \times n_{c}} \rightarrow \mathbb{R}^{n_{y} \times n_{x} \times n_{c}},
\end{equation}
where $n_{y}$ refers to the number of data points in the vertical direction and $n_{x}$ the number of data points in the horizontal direction. The reshape operator is applied to each separate state that occupies the $n_{c}$ input channels. An inverse reshape operator is used to reshape state output data in each output channel into the original vector format
\begin{equation}
    \bm{R}^{-1}: \mathbb{R}^{n_{y} \times n_{x} \times n_{c}} \rightarrow \mathbb{R}^{N \times n_{c}}.
\end{equation}

\section{Expansion coefficient prediction using Gaussian process regression}
Non-intrusive ROMs require a regression model that can accurately predict the expansion coefficients $\bm{a}$ of unrealized solutions given their design parameters $\bm{\mu}$. The regression model is created in the offline stage utilizing the training data. In particular, we use a regression model to approximate a mapping $f(\bm{\mu})$ that outputs the expansion coefficients. A commonly used regression model in non-intrusive ROMs is Gaussian process regression (GPR), a supervised learning method used for predictions of continuous outputs. GPR is also referred to as \emph{Kriging}, and had one of its first uses in the field of geostatistics~\cite{Krige}. The regression model is constructed using training data composed of inputs $\bm{\bar{z}} = [\bm{z}^{1}, \bm{z}^{2} \cdots \bm{z}^{n}]$ and outputs $\bm{\bar{y}} = [y^{1}, y^{2} \cdots y^{n}]$, where each input $\bm{z}^{i} \in \mathcal{P} \subset  \mathbb{R}^{p}$ belonging to an input domain $\mathcal{P}$ corresponds to a single output $y^{i} \in \mathbb{R}$. GPR infers a probability distribution over functions conditioned on the training data which is used for predictions at new inputs. A brief introduction to GPR is given in Section~\ref{sssec:gpr}, and the work of Rasmussen et al.~\cite{rasmussen2003gaussian} can be referred to for a more complete overview.

In ROMs, an individual regression model is used for each coefficient in $\bm{a}$, leading to $k$ different regression models,
\begin{equation}
f_{i}(\bm{\mu}): \mathbb{R}^{p} \rightarrow \mathbb{R}, i \in [1,2, \cdots k]. 
\end{equation}
As a result, ROMs using GPR tend to allow the number of expansion coefficients to be large as long as it does not degrade the quality of the trial manifold. GPR provides reasonable accuracy, is computationally inexpensive, does not require many training samples, and is easy to implement, making it a popular choice in non-intrusive ROMs. While the use of neural networks in non-intrusive ROMs has become more widespread \cite{HESTHAVEN201855,JacquierNN}, finding a sufficient neural network architecture to use for regression is a non-trivial task. Even though using neural networks may offer better performance, in this work our goal is to highlight the advantages of using a nonlinear trial manifold compared to a linear reduced basis. As GPR offers more flexibility, we choose to use it as our regression model.

\subsection{Gaussian process regression}\label{sssec:gpr}
A Gaussian process (GP) is a set of random variables, of which any finite number follow a joint Gaussian distribution. In GPR, it is assumed that data are generated according to a GP with mean function $m$ and covariance function $\kappa$,
\begin{equation}
f(\bm{z}) \sim \text{GP}\left(m(\bm{z}),\kappa(\bm{z},\bm{z}^{*})\right),
\end{equation}
with some added Gaussian noise $\epsilon \sim \mathcal{N}(0,\sigma^{2}_{y})$,
\begin{equation}
    y = f(\bm{z}) + \epsilon.
\end{equation}
Using a finite number of training data $\{\bm{\bar{z}}, \bm{\bar{y}}  \}$, a prior joint Gaussian on the the data and predictions at points $\bm{z}^{*}$ are given by
\begin{equation}
    \begin{bmatrix} \bar{\bm{y}} \\ f(\bm{z}^{*}) \end{bmatrix} \sim \mathcal{N}\left(\begin{bmatrix} m(\bar{\bm{z}}) \\ m(\bar{\bm{z}}) \end{bmatrix}, \begin{bmatrix} \kappa(\bar{\bm{z}}, \bar{\bm{z}}) + \sigma^{2}I & \kappa(\bar{\bm{z}}, \bm{z}^{*}) \\ \kappa(\bm{z}^{*}, \bar{\bm{z}}) & \kappa(\bm{z}^{*}, \bm{z}^{*}) \end{bmatrix} \right).
\end{equation}
Using the properties of conditional Gaussian distributions, the conditional expectation of $f(\bm{z}^{*})$ is given as 
\begin{equation}
    \mathbb{E}(f(\bm{z}^{*})|\bar{\bm{y}}) = \kappa(\bm{z}^{*}, \bar{\bm{z}})(\kappa(\bar{\bm{z}}, \bar{\bm{z}}) + \sigma^{2}_{y}I)^{-1}\left(\bar{\bm{y}}-m(\bar{\bm{z}})\right),
\end{equation}
Where $I$ is the identity matrix. In practice, the mean function $m$ is set to the mean of the training outputs,
\begin{equation}
    m(\bm{\bar{z}}) = \frac{\sum_{i=1}^{n} y^{i}}{n}
\end{equation}
and the inputs are scaled before training to obtain their standard score $\bm{\mathcal{Z}}$ 
\begin{equation}
    \mathcal{Z}^{i}_{j} = \frac{z^{i}_{j}-m(\bm{z}_{j})}{\sigma_{j}},
\end{equation}
where $i$ and $j$ refer to indices of the observation and input entry respectively. There are many kernels that can be chosen for the covariance function. A very common one is the radial basis function (RBF) kernel
\begin{equation}
\kappa(\bm{z},\bm{z}^{*}) = \text{exp}\left(-\frac{d(\bm{z},\bm{z}^{*})^{2}}{2l^{2}}\right).
\end{equation}
Another choice of kernel, and the one that will be used in this work, is the Matern kernel,
\begin{equation}
\kappa(\bm{z},\bm{z}^{*}) = \dfrac{1}{\Gamma(\nu)2^{\nu-1}} \left(\dfrac{\sqrt{2\nu}}{l}d(\bm{z},\bm{z}') \right)^\nu K_{\nu} 
\end{equation}
where $d$ is the Euclidean distance function, $\Gamma$ is the gamma function, and $K_{\nu}$ is the modified Bessel function of the second kind. The set of hyperparameters $\theta$ of the Matern kernel are $l$ and $\nu$, which control the length scale and smoothness respectively. The predictive performance of the regression model is sensitive to the values of the hyperparameters. Gradient-based optimizers are often used to maximize the marginal log-likelihood of the training data to obtain an optimal set of hyperparameters $\theta_{\text{opt}}$
\begin{equation}
    \theta_{\text{opt}} = \underset{\theta}{\text{argmax}} \; \text{log}\,p(\bar{\bm{y}}|\bar{\bm{z}},\theta) = -\frac{1}{2}\bar{\bm{y}}^{T}(\kappa(\bar{\bm{z}}, \bar{\bm{z}}) + \sigma^{2}I)^{-1} -\frac{1}{2}\text{log}\,|\kappa(\bar{\bm{z}}, \bar{\bm{z}}) + \sigma^{2}I| - \frac{n}{2}\text{log}\,2\pi .
\end{equation}

\section{Offline and online stages}
This section describes the offline training and online evaluation stages of both the POD and CAE based ROMs with GPR as a regression model. The combined models consisting of both the offline and online stages are referred to respectively as POD-GPR and CAE-GPR. The offline stage is run first and is computationally expensive, while the online stage allows for rapid prediction of full-order models. Both models share a step of obtaining full-order snapshots of solutions evaluated at a set of design parameters $\bm{\mathcal{U}}_{\text{train}}$ and assembling them into a snapshot matrix $\bm{S}$. The offline stage of the POD-GPR method involves calculating a truncated SVD of the snapshot matrix to find the POD basis $\bm{\Psi}$. The set of expansion coefficients $\bm{A}_{\text{train}}$ is obtained for the training set using the POD basis and $k$ GPR models $\mathcal{\bm{F}} = [f_{1}(\bm{\mu}), f_{2}(\bm{\mu}), \cdots f_{k}(\bm{\mu})]$ trained on $\bm{\mathcal{T}} = \{\bm{\mathcal{U}}_{\text{train}},\bm{A}_{\text{train}}\}$. After training, the GPR models are saved for use in the online stage, where unseen parameters $\bm{\mu}^{*}$ are evaluated to approximate the expansion coefficients $\bm{\tilde{a}}^{*}$. Matrix-vector multiplication of the POD-basis and expansion coefficients is then used to obtain an approximate solution $\bm{\tilde{x}}$. The POD-GPR method is outlined in Algorithm~\ref{alg:pod-gpr}.

\begin{algorithm}[!htbp]
\caption{Offline and online stages of POD-GPR method}\label{alg:pod-gpr}
\begin{algorithmic}[1]
\Function{PODGPR\_ OFFLINE}{$\bm{\mathcal{U}}_{\text{train}}$}
\State Compute high-fidelity solutions for $\bm{\mu} \in \bm{\mathcal{U}_{\text{train}}}$ by solving FOM and assemble into $\bm{S}$
\State Calculate truncated SVD of snapshot matrix to obtain POD basis $\bm{\Psi}$
\State Calculate expansion coefficients for training data $\bm{A}_{\text{train}} = \left( \bm{\Psi}^{T}\bm{S}\right)^{T}$
\State Train $k$ GPR models $\mathcal{\bm{F}} = [f_{1}(\bm{\mu}), f_{2}(\bm{\mu}), \cdots f_{k}(\bm{\mu})]$ for each expansion coefficient in $ \{\bm{\mathcal{U}}_{\text{train}},\bm{A}_{\text{train}}\}$ 
\State \Return $\left(\bm{\Psi}, \mathcal{F} \right)$
\EndFunction
\end{algorithmic}
\[\]
\begin{algorithmic}[1]
\Function{PODGPR\_ ONLINE}{$\bm{\mu}^{*}, \bm{\Psi}, \mathcal{F}$}
\State Evaluate expansion coefficients $\bm{\tilde{a}}^{*} = \mathcal{F}(\bm{\mu}^{*})$
\State Predict full-order solution $\bm{\tilde{x}}^{*} = \bm{\Psi}\bm{\tilde{a}}^{*}$
\State \Return $\bm{\tilde{x}}^{*}$
\EndFunction
\end{algorithmic}
\end{algorithm}

\begin{algorithm}[!htbp]
\caption{Offline and online stages of CAE-GPR method}\label{alg:cae-gpr}
\begin{algorithmic}[1]
\Function{CAEGPR\_ OFFLINE}{$\bm{\mathcal{U}}_{\text{train}}, \bm{\mathcal{U}}_{\text{\text{val}}},  \bm{\mathcal{C}}, \bm{R} $}
\State Compute high-fidelity solutions for $\bm{\mu} \in \bm{\mathcal{U}_{\text{train}}},  \bm{\mathcal{U}_{\text{\text{val}}}} $ by solving FOM and assemble into $\bm{S}_{\text{train}}, \bm{S}_{\text{\text{val}}} $
\State Apply reshape operator $\bm{R}$ to $\bm{S}_{\text{train}}, \bm{S}_{\text{\text{val}}} $ to obtain $\bm{X}_{\text{train}}, \bm{X}_{\text{\text{val}}} $ 
\State Train convolutional autoencoder with architecture $\bm{\mathcal{C}}$ on $\{\bm{X}_{\text{train}}, \bm{X}_{\text{train}}\}$ while monitoring loss on $\{\bm{X}_{\text{\text{val}}}, \bm{X}_{\text{\text{val}}}\}$
\State Calculate expansion coefficients for training data $\bm{A}_{\text{train}} = g_{\text{enc}}\left(\bm{X}_{\text{train}} \right)$
\State Train $k$ GPR models $\mathcal{\bm{F}} = [f_{1}(\bm{\mu}), f_{2}(\bm{\mu}), \cdots f_{k}(\bm{\mu})]$ for each expansion coefficient in $ \{\bm{\mathcal{U}}_{\text{train}},\bm{A}_{\text{train}}\}$ 
\State \Return $\left(g_{\text{dec}}, \mathcal{F} \right)$
\EndFunction
\end{algorithmic}
\[\]
\begin{algorithmic}[1]
\Function{CAEGPR\_ ONLINE}{$\bm{\mu}^{*},g_{\text{dec}}, \mathcal{F}, \bm{R}^{-1}$}
\State Evaluate expansion coefficients $\bm{\tilde{a}}^{*} = \mathcal{F}(\bm{\mu}^{*})$
\State Predict full-order solution $\bm{\tilde{X}}^{*} = g_{\text{dec}}\left(\bm{\tilde{a}}^{*}\right)$
\State Apply inverse reshape operator $\bm{R}^{-1}$ to  $\bm{\tilde{X}}^{*}$ to obtain  $\bm{\tilde{x}}^{*}$
\State \Return $\bm{\tilde{x}}^{*}$
\EndFunction
\end{algorithmic}
\end{algorithm}
The CAE-GPR method also involves evaluating additional full-order solutions at a set of validation design parameters $\bm{\mathcal{U}}_{\text{\text{val}}}$ which are used to monitor the validation loss during training so early stopping can be implemented as a regularization method. A convolutional autoencoder architecture $\bm{\mathcal{C}}$ is also required and the decoder $g_{\text{dec}}$ is saved for use in the online stage. Approximate expansion coefficients $\bm{\tilde{a}}$ are passed to the decoder to obtain approximate solutions $\bm{\tilde{X}}.$ The reshape operator $\bm{R}$ is also required to make the training and validation data compatible with the architecture of the CAE in the offline stage, while the inverse reshape operator $\bm{R}^{-1}$ is needed in the online stage to reshape approximated solutions to their original format. The CAE-GPR method is outlined in Algorithm~\ref{alg:cae-gpr}. Training the convolutional autoencoder makes the offline stage of CAE-GPR more expensive than POD-GPR, while the online costs for both models are similar.

\section{Numerical results}
This section compares the performance of the POD-GPR and CAE-GPR methods on a geometrically and physically parameterized lid-driven cavity problem which simulates laminar flow using the steady incompressible Navier-Stokes equations with OpenFOAM~\cite{weller1998tensorial}, an open-source toolbox for multiphysics simulation. The FOM quantities of interest are $u$ and $v$, the components of the velocity in the horizontal and vertical directions respectively. The autoencoder is constructed using TensorFlow~\cite{tensorflow2015-whitepaper}, and scikit-learn~\cite{Pedregosa2011} is used to implement GPR. 500 design parameters are generated using Latin hypercube sampling~\cite{lhs}, a statistical method that aims to maximize the distance and minimize the correlation amongst produced samples. The metric of performance used to compare the ROMs is the relative $l^{2}$ error $\epsilon_{\text{ROM}}$ between the FOM state $\bm{x}$ and the ROM approximated state $\bm{\tilde{x}}$
\begin{equation}
        \epsilon_{\text{ROM}} =  \dfrac{\norm{\bm{x}^{i} - \tilde{\bm{x}}^{i}}^{2}} {\norm{\bm{x}^{i}}^2}.
\end{equation}
Similarly, the relative projection error $\epsilon_{\text{Proj}}$ for both methods is also reported between the FOM state and the projected state $\bm{\hat{x}}$ to assess how accurately the expansion coefficients are interpolated as well as to provide a lower bound for the ROM prediction errors.

\begin{equation}
        \epsilon_{\text{Proj}} =  \dfrac{\norm{\bm{x}^{i} - \hat{\bm{x}}^{i}}^{2}} {\norm{\bm{x}^{i}}^2}.
\end{equation}

A five-fold cross-validation approach is used to assess the performance of the ROM over the entire dataset, creating five folds of the dataset containing 400 training samples and 100 testing/validation samples. These 100 samples are split evenly into 50 testing and 50 validation samples, which are used to monitor the autoencoder loss during training. While every sample in the dataset is used for training, only half are used for prediction. An average cross-validation error is reported for both the ROM prediction and projection errors over all of the prediction points. The validation samples are not used for the POD-GPR ROM, which uses an individual ROM for both $u$ and $v$.

The CAE that is used has two input channels, one for each of the velocity components. The expansion coefficients are used to approximate both $u$ and $v$.  The encoder conists of a combination of convolutional, pooling, and fully connected layers. The decoder consists of fully connected and transpose-convolutional layers. The full details of the network architecture can be found in the appendix. Min-max scaling is used independently on $u$ and $v$ before training and the CAE outputs are then scaled back to their original range after prediction. A maximum number of 7500 training epochs are used, and early stopping is enforced if the validation loss fails to decrease over 500 epochs. A mini-batch size of $b = 8$ is used for training and the mean squared error loss function is used. The Adam optimizer is used with an initial learning rate of of $\eta = 3 \times 10^{-4}$. Each layer has its weights initialized with the He normal initializer. All of the layers with the exception of the output use a leaky ReLU activation function 
with $\alpha$ = 0.25. The output layer uses the sigmoid activation function, which scales into the range [0,1], ensuring that the outputs can be scaled back to their original range for prediction
\begin{equation}
\phi(x) = \frac{1}{1+e^{-x}}.
\end{equation}

As the FOM is inexpensive to solve, offline computational costs related to simulation and POD-GPR are not reported as well as all online costs. Training and validation losses against the number of epochs at selected folds of the data for different values of $k$ are presented in the appendix in addition to computational costs.

\subsection{Steady incompressible Navier-Stokes equations}
Steady incompressible laminar flow is simulated using simpleFoam, a standard OpenFOAM solver, by solving the Navier-Stokes equations,
\begin{equation}
\label{eqn_mass}
\int_S \overrightarrow{U} \cdot \dif \overrightarrow{S}=0,
\end{equation}
\begin{equation}
  \label{eqn_momentum}
\int_S  \overrightarrow{U} \overrightarrow{U} \cdot \dif \overrightarrow{S}  + \int_V \nabla p  \dif V -\nu \int_S (\nabla \overrightarrow{U}+\nabla \overrightarrow{U}^T) \cdot \dif \overrightarrow{S}  = 0,
\end{equation}
where $\overrightarrow{U}$ = [$u$, $v$] is the velocity vector and $u$ and $v$ are the velocity components in the $x$ and $y$ directions respectively, $\overrightarrow{S}$ is the face-area vector, $V$ is the volume; $\nu$ is the kinematic viscosity, and $p$ is the pressure.
The continuity and momentum equations are discretized over the computational domain by using the finite-volume method (FVM). Both equations are coupled through the semi-implicit method for pressure-linked equations (SIMPLE) algorithm~\cite{Patankar1972} along with Rhie--Chow interpolation~\cite{Rhie1983}. The SIMPLE algorithm is iteratively repeated until a residual tolerance of $1 \times 10^{-6}$ is reached for both $\overrightarrow{U}$ and $p$.

\subsection{Lid-driven cavity problem}
\begin{figure}[!t]
\centering
\subfigure[Lid-driven cavity boundary conditions.]{\includegraphics[width=0.49\textwidth]{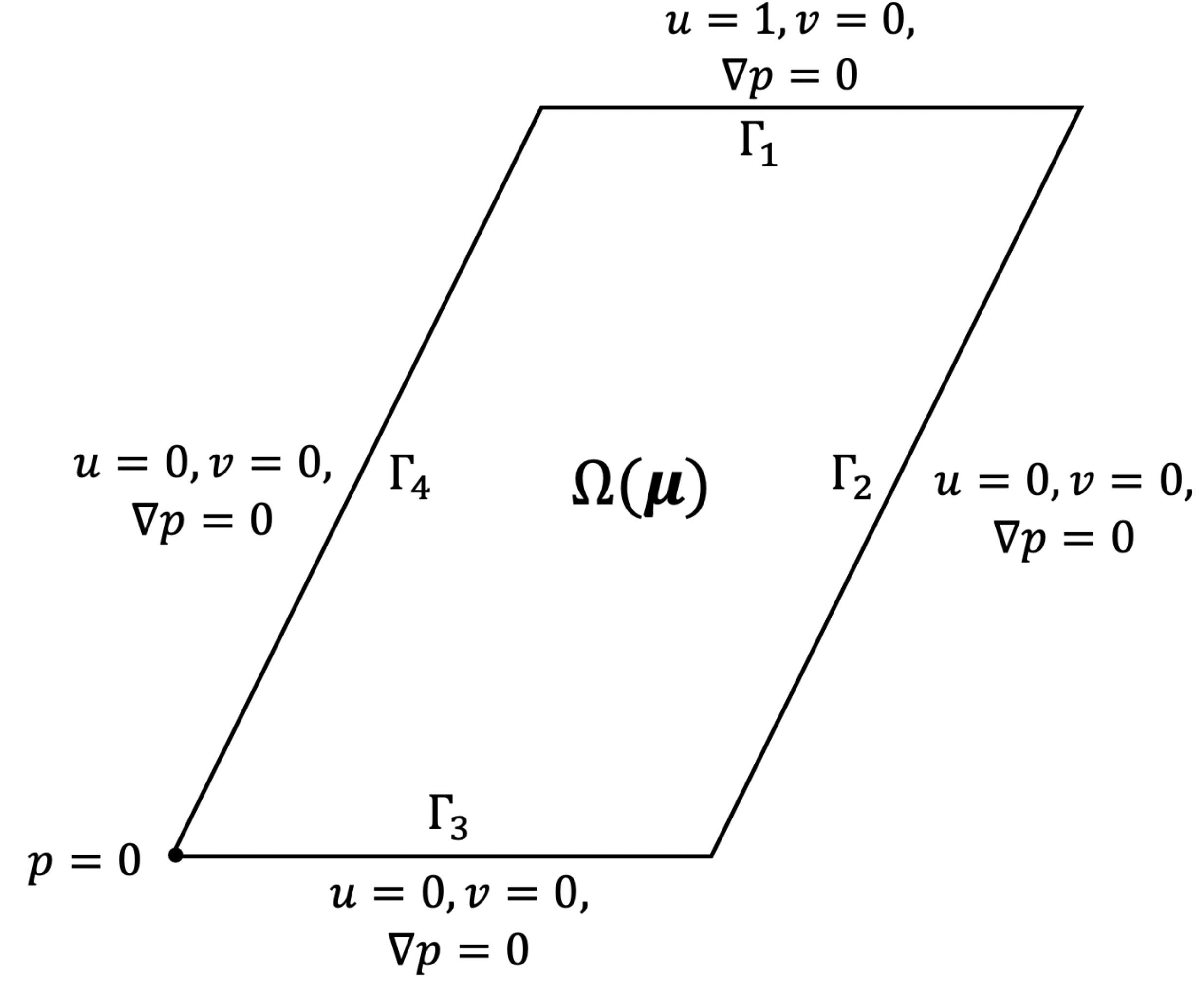}}
\subfigure[Lid-driven cavity design parameters.]{\includegraphics[width=0.49\textwidth]{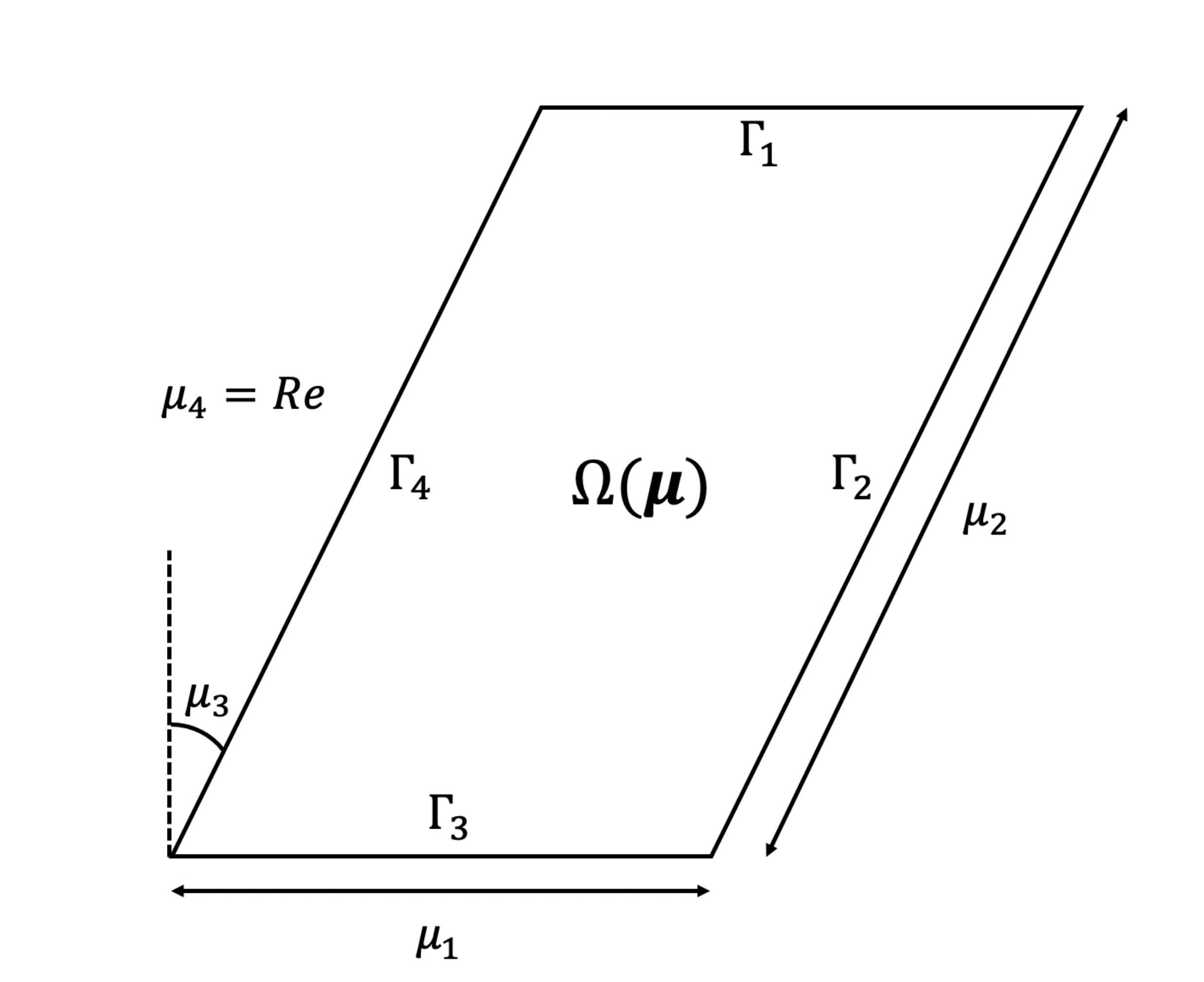}}
\caption{Schematics describing the lid-driven cavity problem.}
\label{fig:cavity_schematic}
\end{figure}

\begin{figure}[!t]
\centering
\subfigure[$\bm{\mu} = (1.5, \frac{3}{\sqrt{3}}, -\frac{\pi}{6}, 200)$]{\includegraphics[scale=.125]{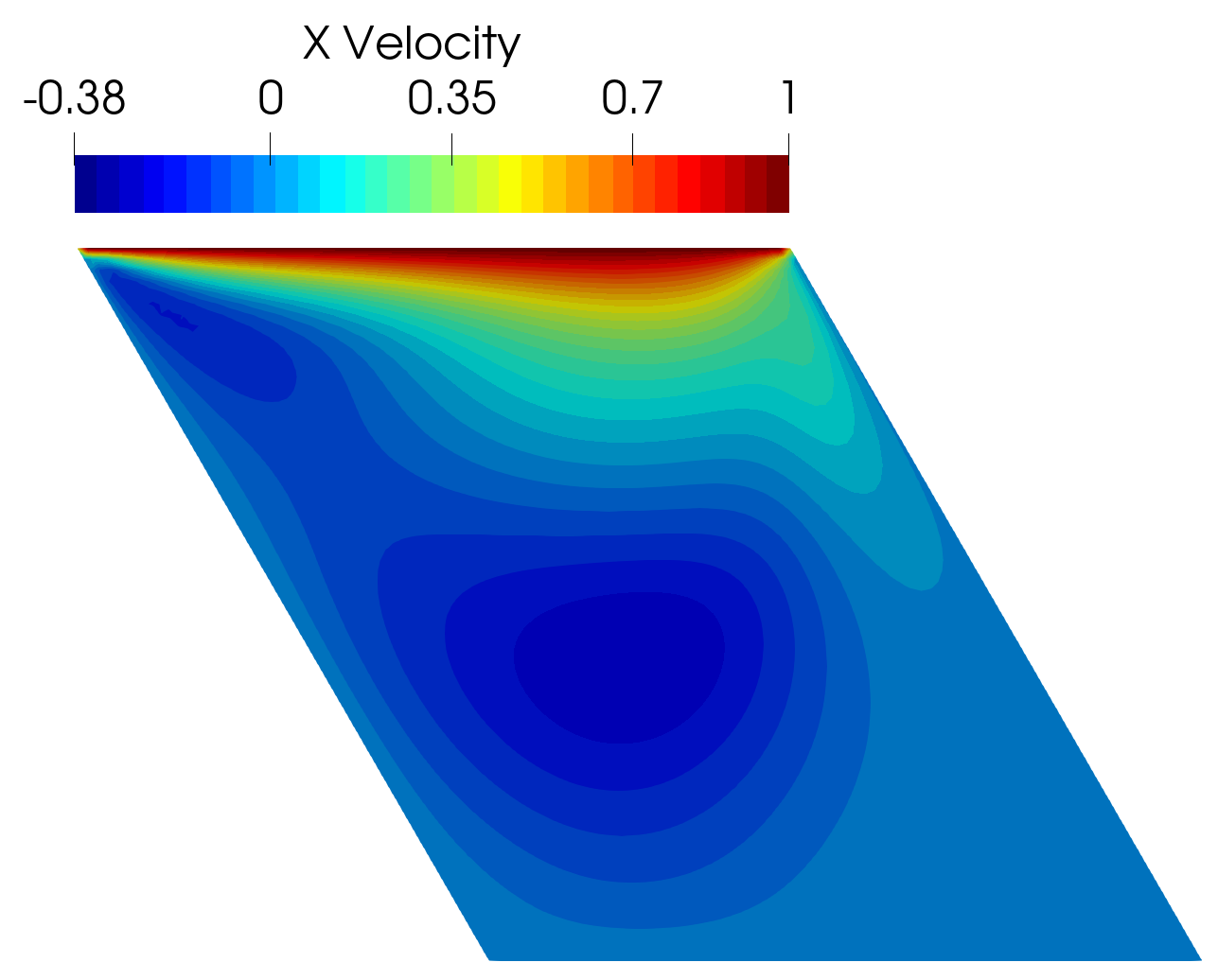}}
\hspace{0.2cm}
\subfigure[$\bm{\mu} = (1.5, 1.5, 0, 300)$]{\includegraphics[scale=.125]{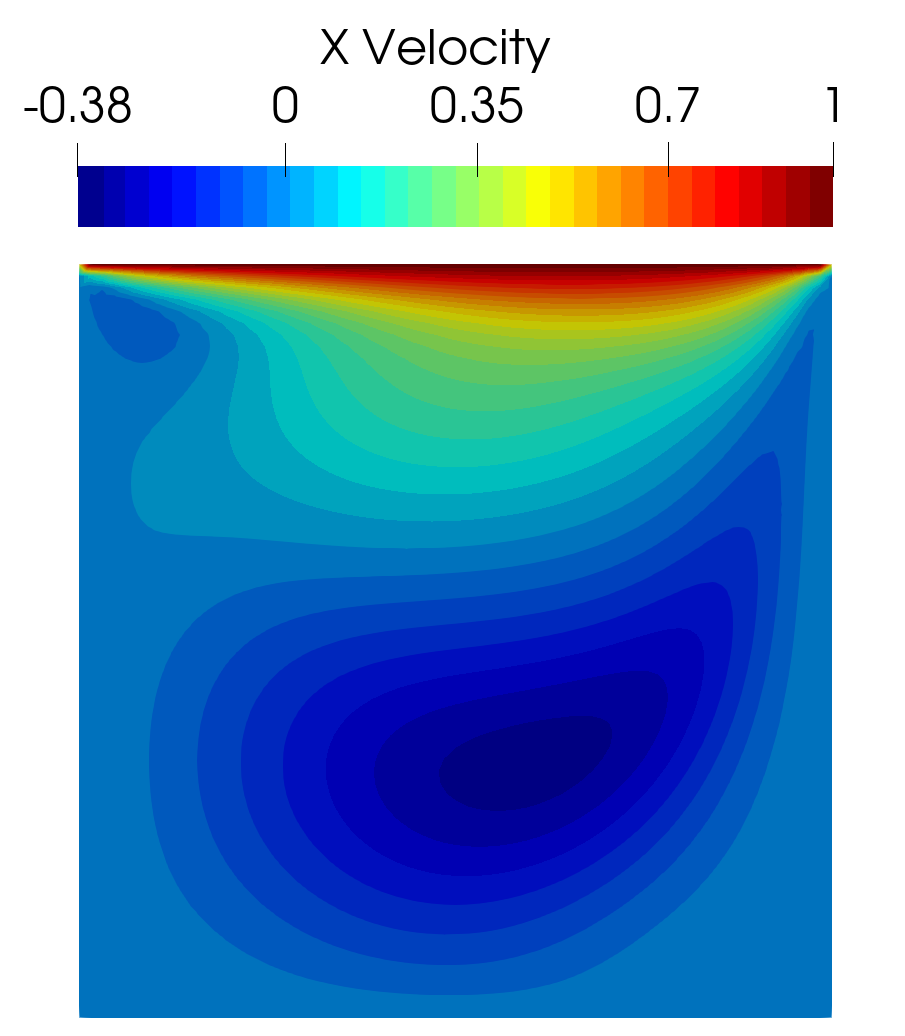}} 
\hspace{0.2cm}
\subfigure[$\bm{\mu} = (1.5, \frac{3}{\sqrt{3}}, \frac{\pi}{6}, 400)$]{\includegraphics[scale=.125]{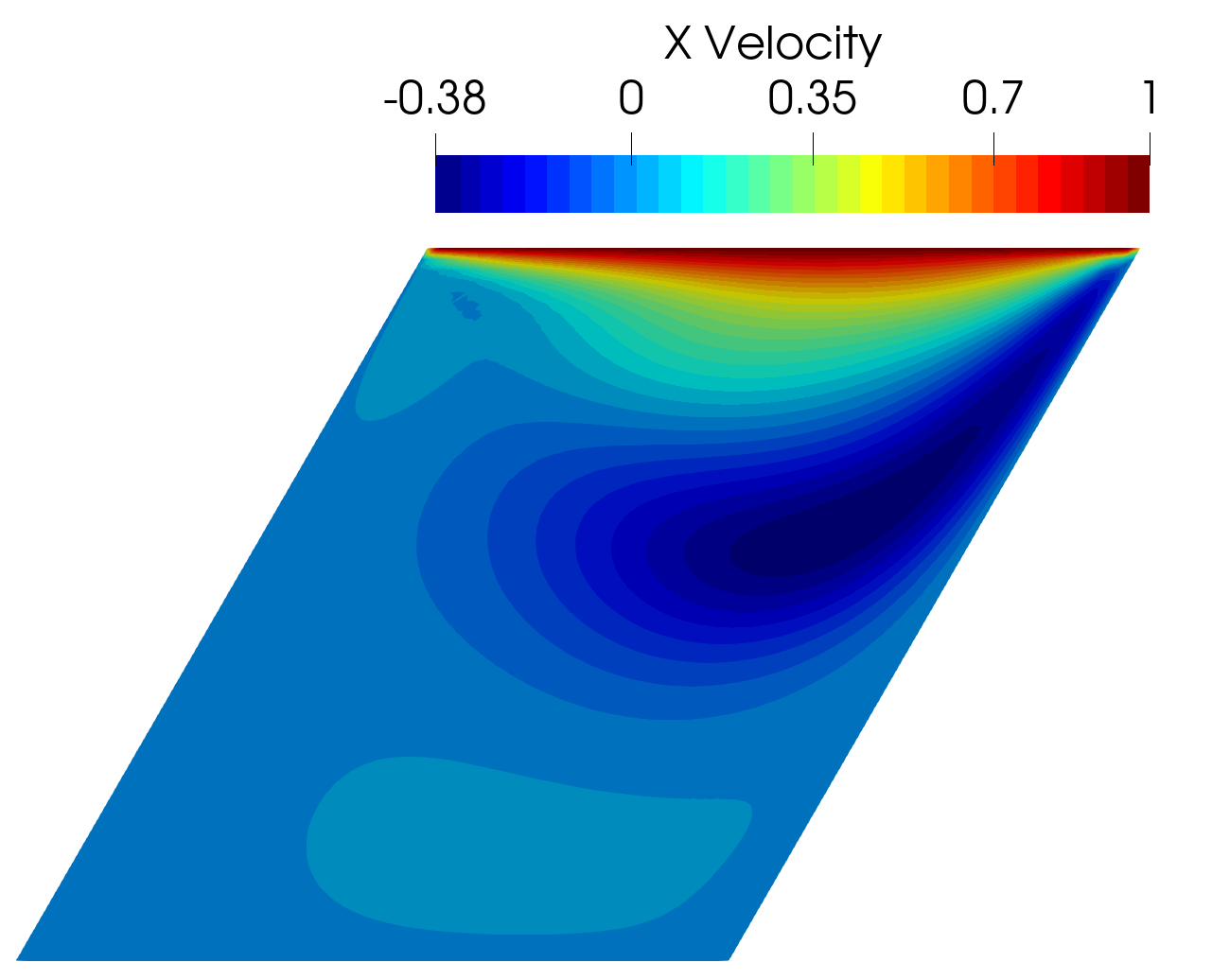}}
\subfigure[$\bm{\mu} = (1.5, \frac{3}{\sqrt{3}}, -\frac{\pi}{6}, 200)$]{\includegraphics[scale=.125]{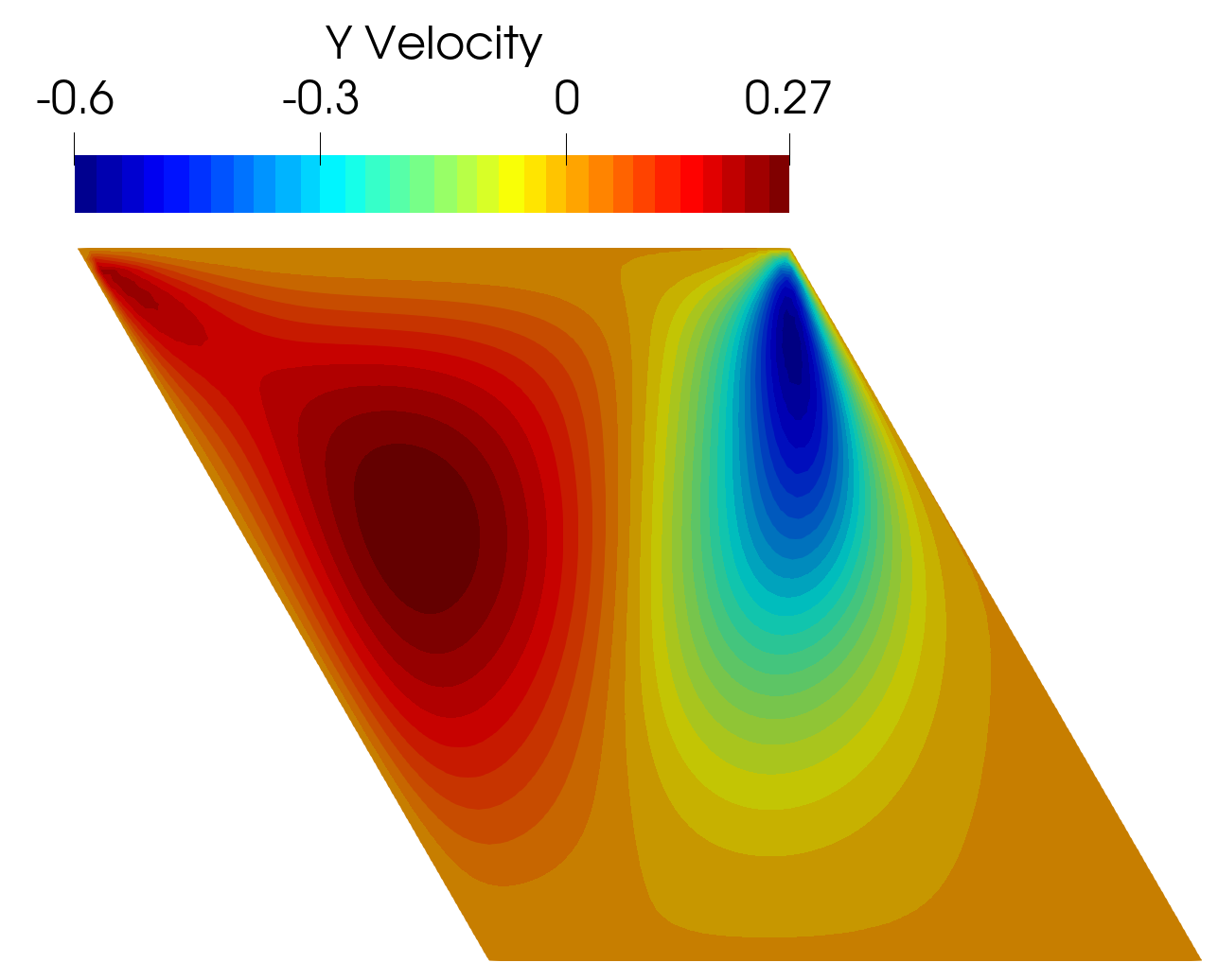}} 
\hspace{0.2cm}
\subfigure[$\bm{\mu} = (1.5, 1.5, 0, 300)$]{\includegraphics[scale=.125]{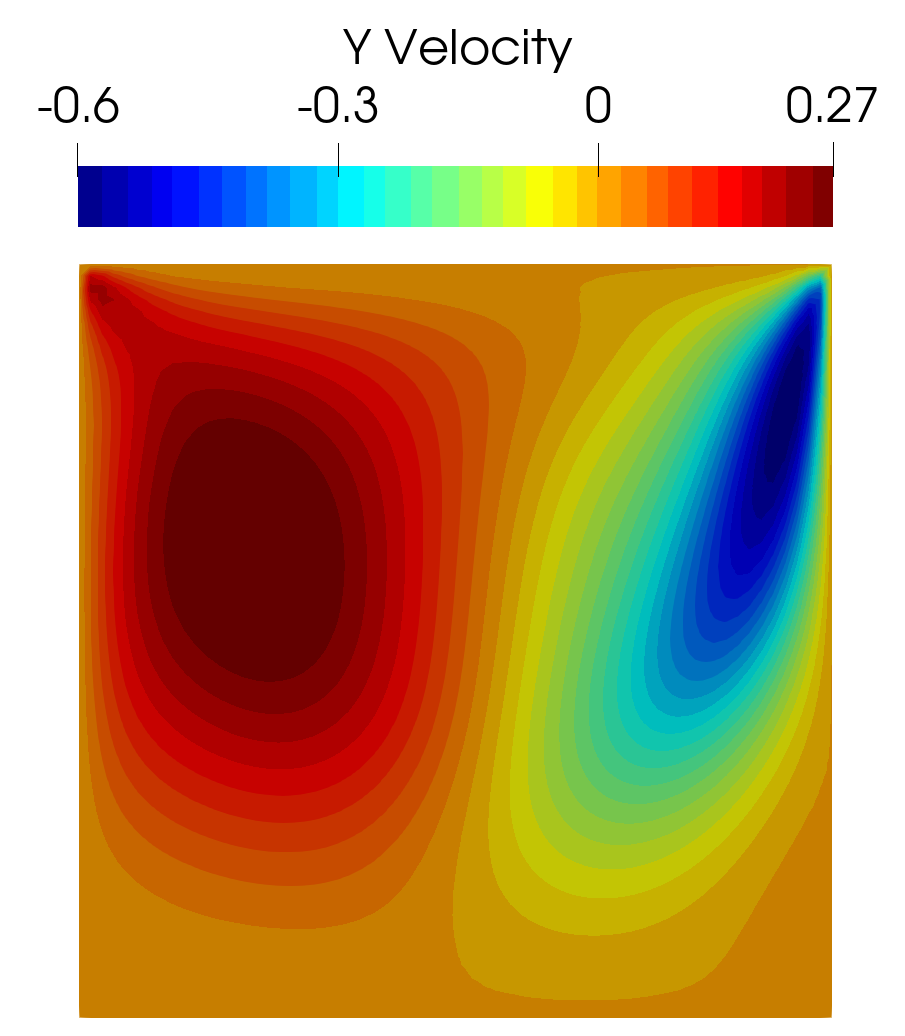}} 
\hspace{0.2cm}
\subfigure[$\bm{\mu} = (1.5, \frac{3}{\sqrt{3}}, \frac{\pi}{6}, 400)$]{\includegraphics[scale=.125]{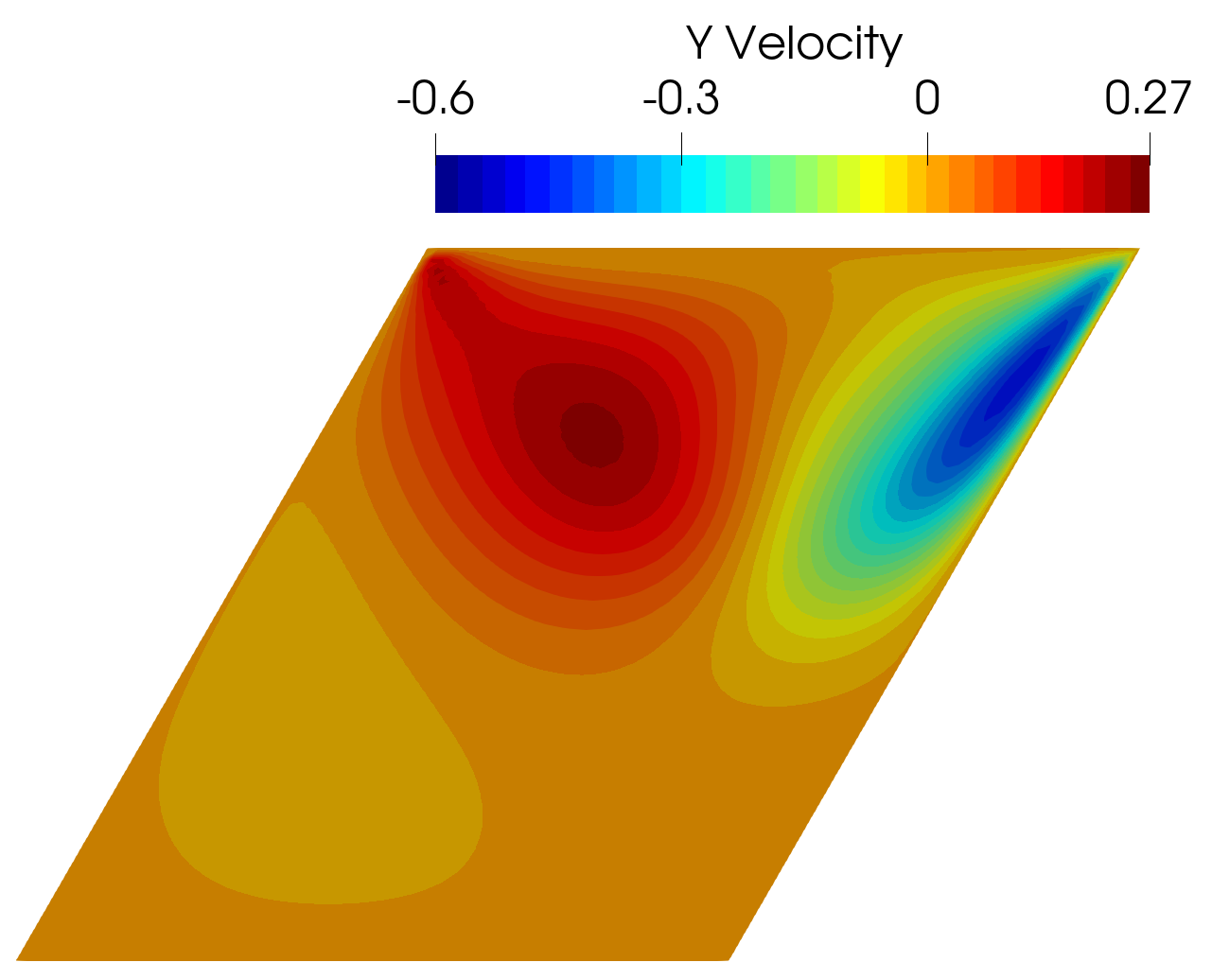}}
\caption{Contours of $u$ (top) and $v$ (bottom) for the lid-driven cavity problem at three different sets of design parameters.}
\label{fig:cavity_contours}
\end{figure}

The numerical example used in this work is a physically and geometrically parameterized lid-driven cavity flow, a popular benchmark problem for CFD solvers. Three parameters control the computational domain $\Omega$ and one parameter controls the kinematic viscosity through the Reynolds number. A version of this problem has previously appeared in a work by Hesthaven and Ubbiali~\cite{HESTHAVEN201855}. Figure~\ref{fig:cavity_schematic} shows the boundary conditions on each edge $\Gamma_{i}, i \in [1,2,3,4]$ of the domain; $u$, $v = 0$ on all of the edges except $\Gamma_{1}$, where $u = 1, v = 0$. The pressure gradient, $\nabla{p},$ is set to 0 on all of the edges. The reference pressure is set to 0 on the bottom left corner of the domain. The parameterization of the geometry is also shown, involving three parameters which change the length of the horizontal $(\mu_{1})$ and slanting edges $(\mu_{2})$ as well as the slanting angle $(\mu_{3})$. The Reynolds number, $Re$ $(\mu_{4})$, is the fourth parameter, and is related to the kinematic viscosity $\nu$ as
\begin{equation}
    Re = \frac{\text{max}(\mu_{1}, \mu_{2})}{\nu(\bm{\mu})}.
\end{equation}
The design parameter combinations are generated using Latin hypercube sampling with the following bounds for each parameter
\begin{align*}
& \mu_{1} \in [1, 2], \\
& \mu_{2} \in [1, 2], \\ 
& \mu_{3} \in [-\frac{\pi}{4}, \frac{\pi}{4}], \\ 
& \mu_{4} \in [100, 600].
\end{align*}
The computational mesh consists of $64 \times 64$ cells uniformly distributed in the $x$ and $y$ directions and one cell spanning the $z$ direction, resulting in $N = 4096$ and a reshape operator $\bm{R}$ with $n_{y}, n_{x}$ = 64 and $n_{c}$ = 2. The full-order states of $u$ and $v$ are used to compare the performance of the POD-GPR and CAE-GPR methods, contours of which are shown in Figure~\ref{fig:cavity_contours} at three different sets of design parameters. A sharp gradient in $u$ exists at the top of the domain, and a vortex moves throughout the cavity as the design parameters change. This vortex is also shown moving throughout the cavity shown in the contours of $v$, varying in shape and size with the design parameters. The relationship between both $u$ and $v$ and $\bm{\mu}$ is shown to be highly nonlinear, making this a difficult prediction problem in the context of ROMs. 

\begin{figure}[!htbp]
\centering
\subfigure{\includegraphics[width=0.49\textwidth]{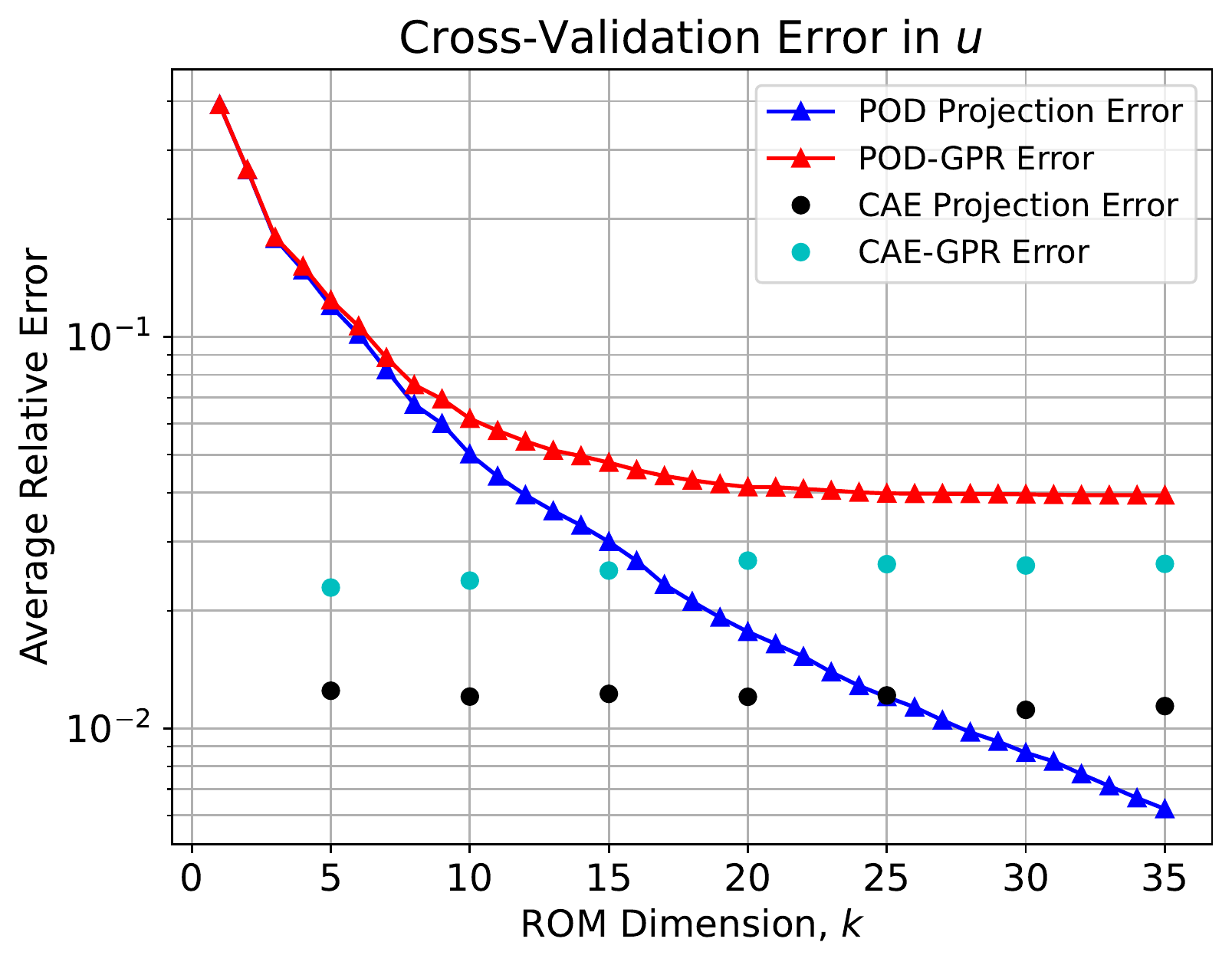}}
\subfigure{\includegraphics[width=0.49\textwidth]{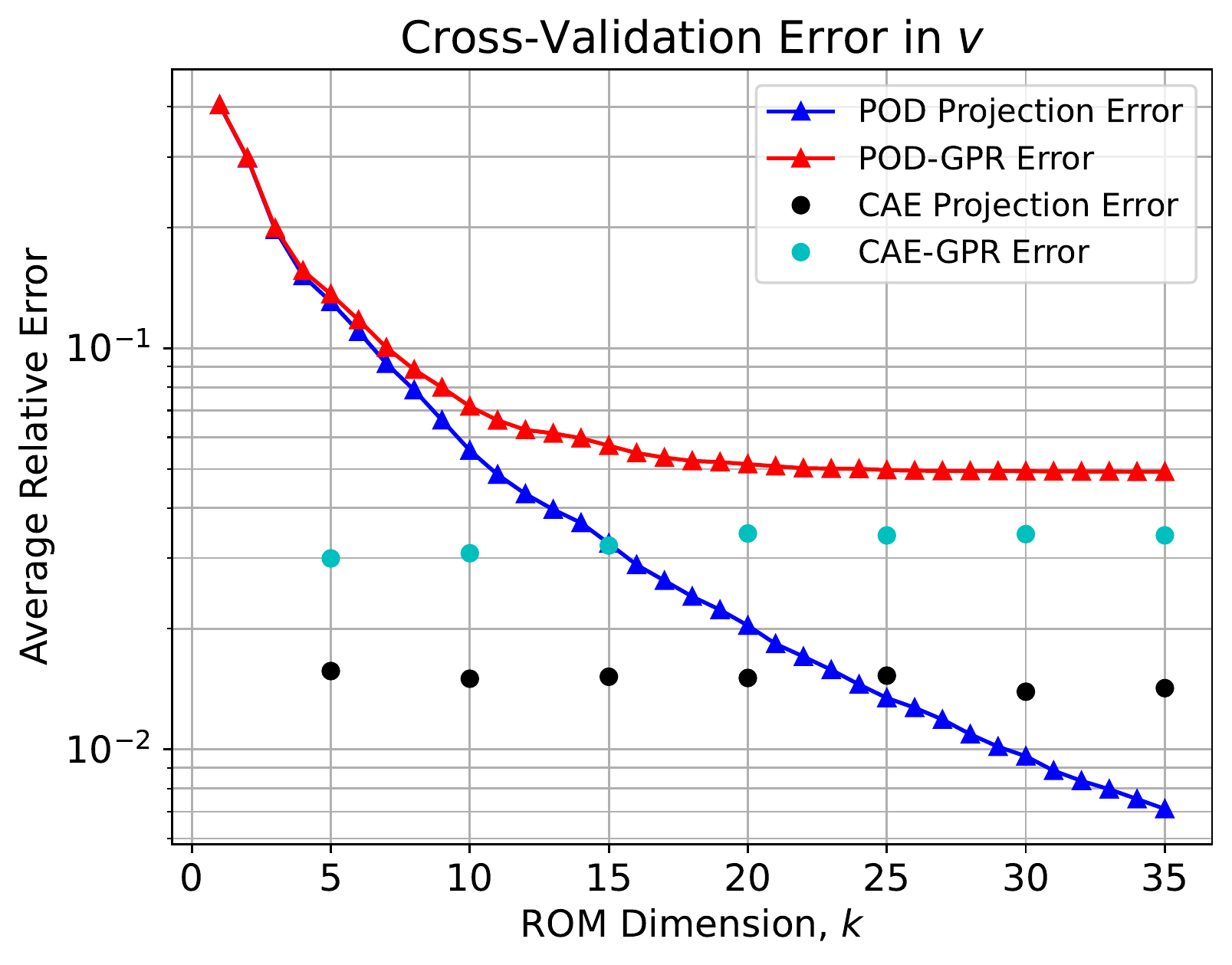}}
\caption{Plots of the cross-validation prediction and projection errors in $u$ and $v$ for both ROMs at different values of $k$.} 
\label{fig:err_cv}
\end{figure}

The POD-GPR ROM is evaluated for both projection and prediction errors at ROM dimensions $k \in [1, 2, \cdots 35]$ while the CAE-GPR ROM is similarly evaluated at $k \in [5, 10, \cdots 35]$. Figure~\ref{fig:err_cv} shows $\bar{e}_{CV}$ for $u$ and $v$; over all evaluated ROM dimensions, the CAE-GPR ROM exhibits higher predictive performance over the data set. The projection error provided by the CAE does not vary much with $k$, while the prediction error tends to increase slightly with $k$. This is an expected result, since the projection errors are very similar, and prediction errors in the individual expansion coefficients will have a cascading effect. The cross-validation projection error from POD continues to decay after $k = 35$, while the POD-GPR prediction error flattens out at around $k = 20$. The cross-validation projection error produced by the CAE for both $u$ and $v$ is lower than that of POD until around $k = 25$; even with a higher projection error, the predictive performance offered by CAE-GPR exceeds that of POD-GPR for $k \geq 25$. The CAE offers a set of expansion coefficients that are more easily interpolated when using GPR compared to POD, which sees its predictive performance stall after $k$ reaches a certain value, a commonly found result for POD-GPR based ROMs~\cite{MrosekROM}. In addition to giving better performance in terms of both projection and prediction for low values of $k$, the use of a nonlinear trial manifold for ROM construction offers a more robust relationship between the design parameters and expansion coefficients. The difference in projection and prediction errors from POD-GPR is almost 0 at low values of $k$, but rapidly increases as $k$ grows, suggesting that the individual expansion coefficients become harder to interpolate as the corresponding singular values decay. 

\begin{figure}[!htbp]
\centering
\subfigure{\includegraphics[width=0.49\textwidth]{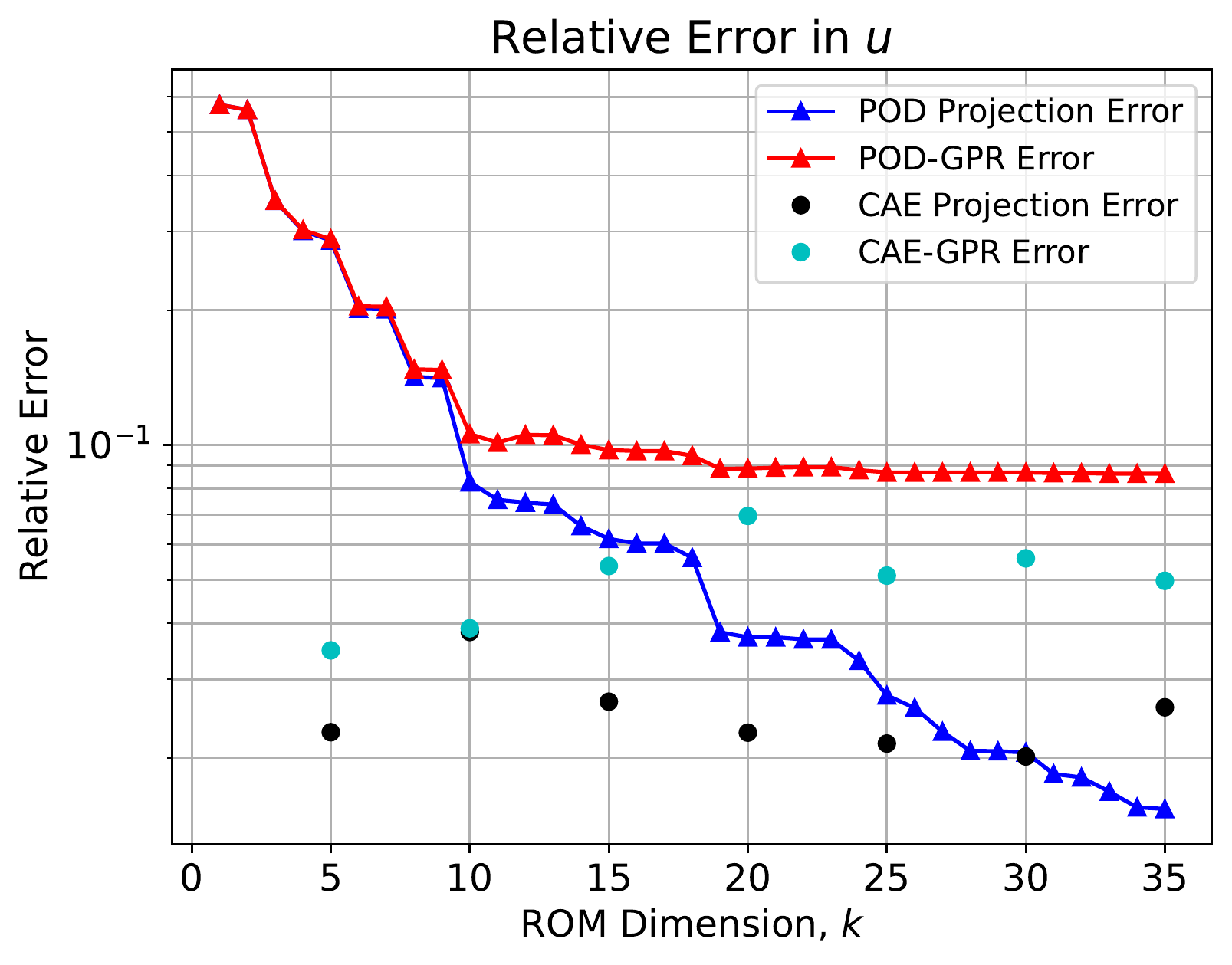}}
\subfigure{\includegraphics[width=0.49\textwidth]{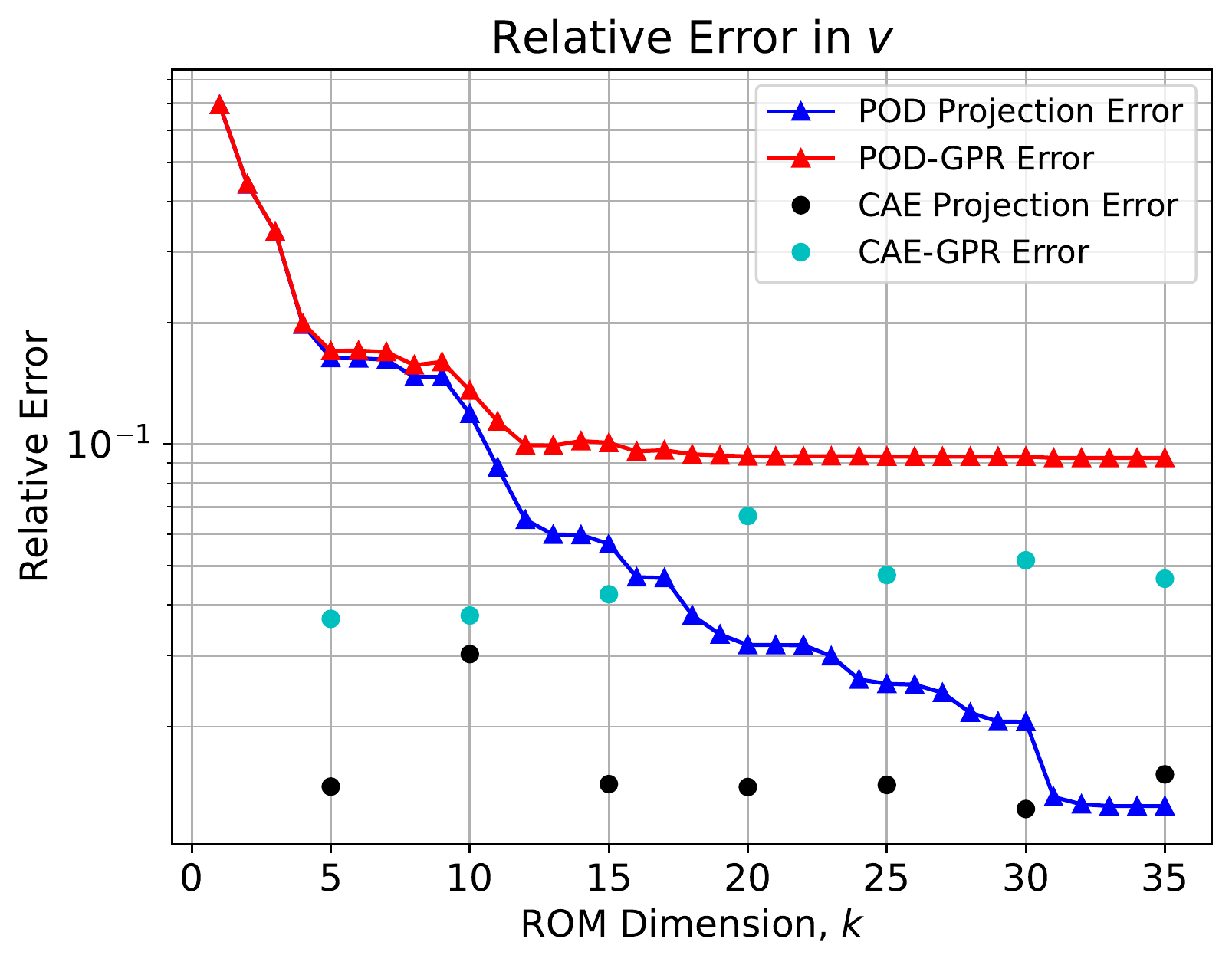}}
\caption{Plots of the prediction and projection errors in $u$ and $v$ for both ROMs at $\bm{\mu} = (1.167, 1.997, -0.4665, 555.5)$ at different values of $k$.} 
\label{fig:err_322}
\end{figure}

\begin{figure}[!t]
\centering
\subfigure[$\bm{u}, \text{Ground Truth}$]{\includegraphics[scale=0.105]{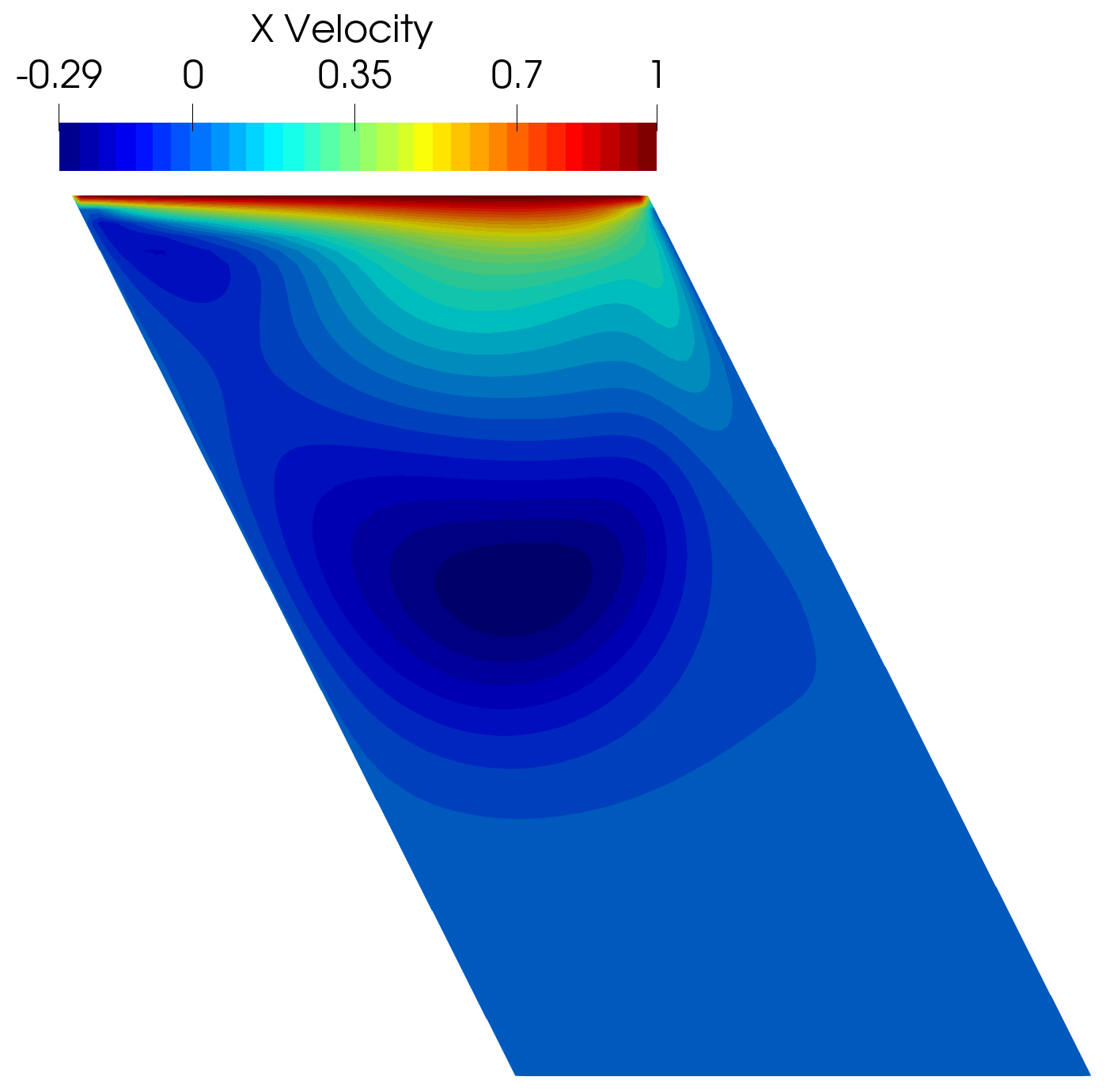}}
\hspace{0.2cm}
\subfigure[$\bm{u}, \text{POD-GPR Difference}$]{\includegraphics[scale=0.105]{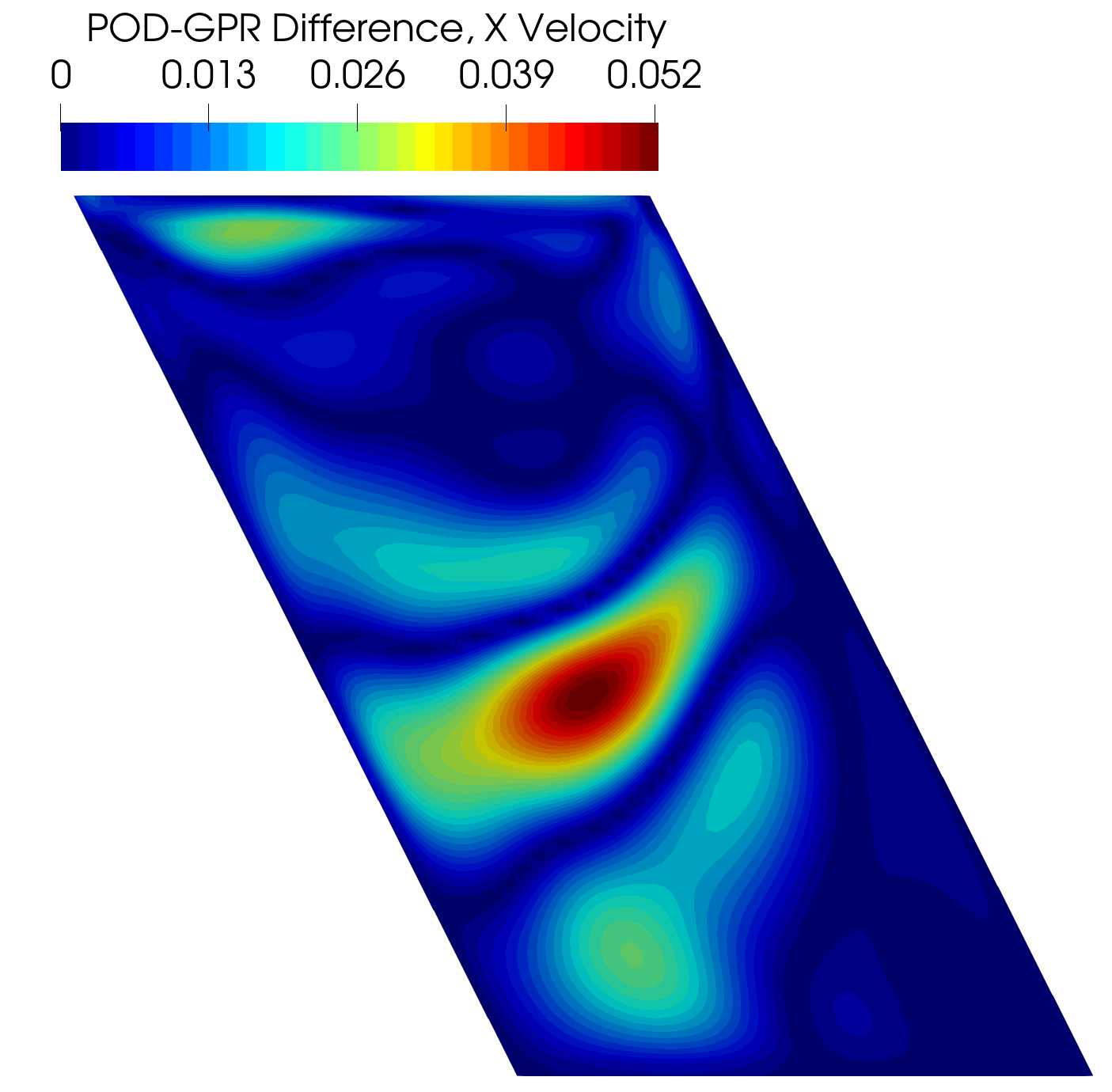}} 
\hspace{0.2cm}
\subfigure[$\bm{u}, \text{CAE-GPR Difference}$]{\includegraphics[scale=0.105]{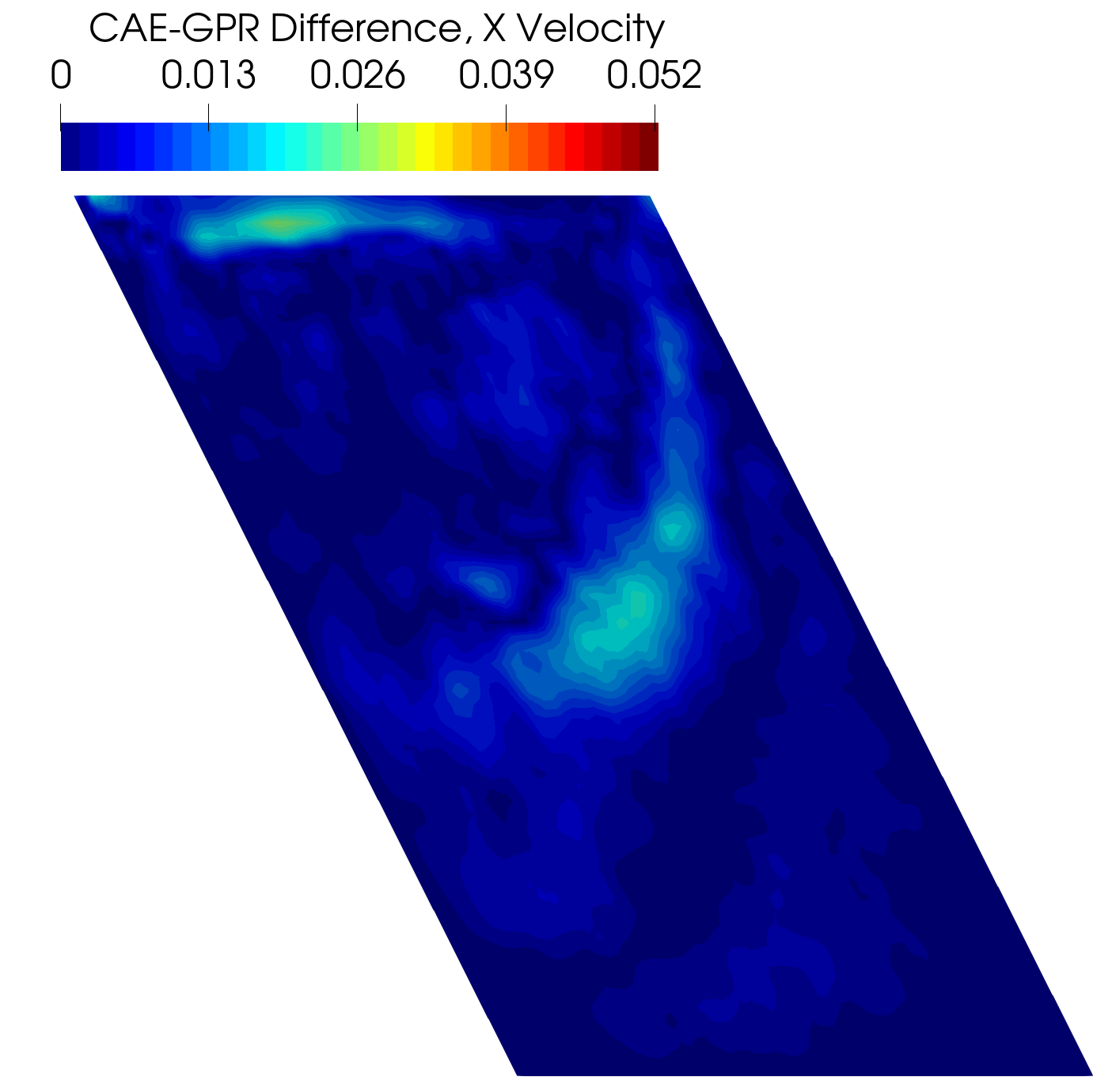}}
\subfigure[$\bm{v}, \text{Ground Truth}$]{\includegraphics[scale=0.105]{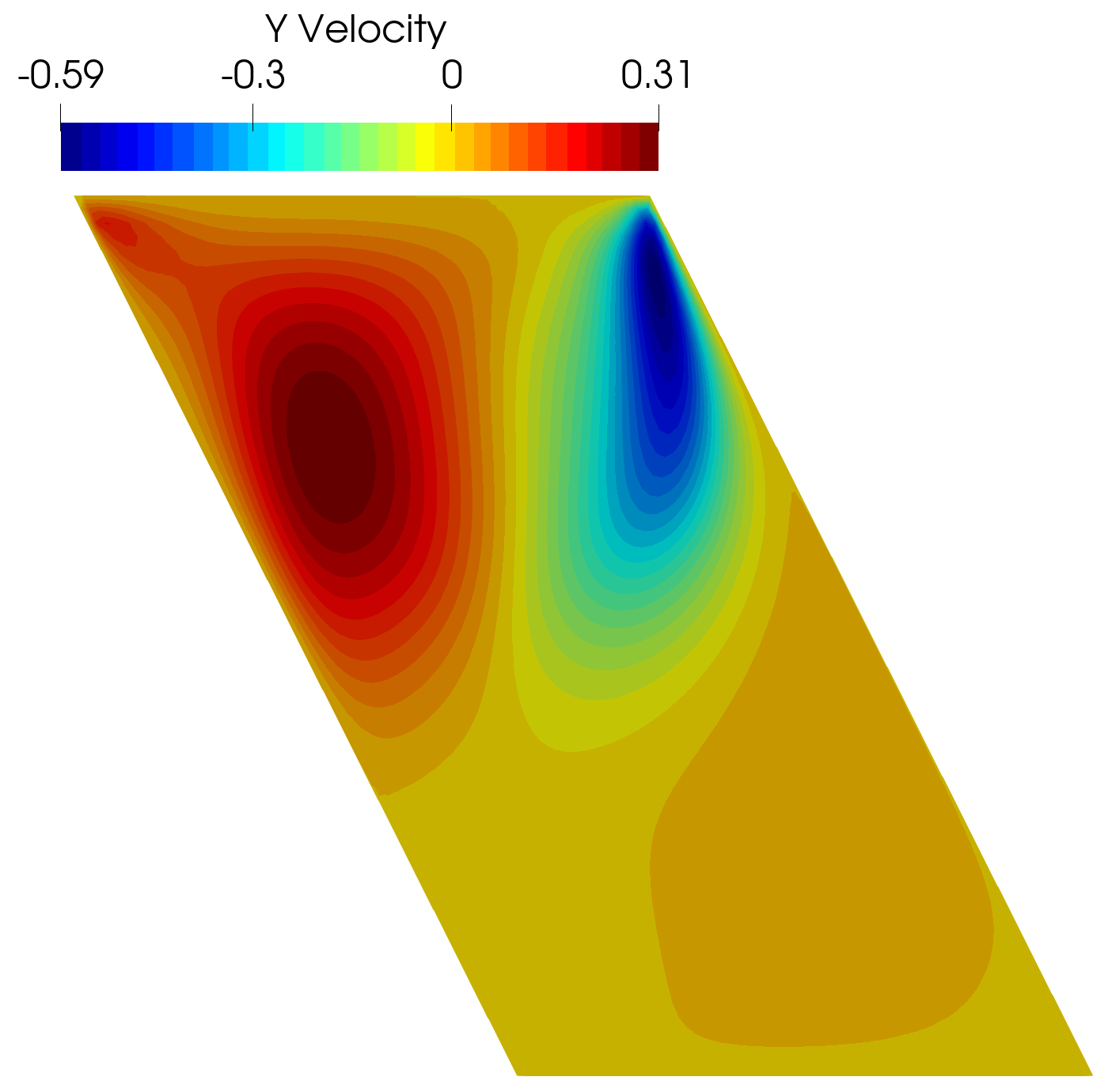}}
\hspace{0.2cm}
\subfigure[$\bm{v}, \text{POD-GPR Difference}$]{\includegraphics[scale=0.105]{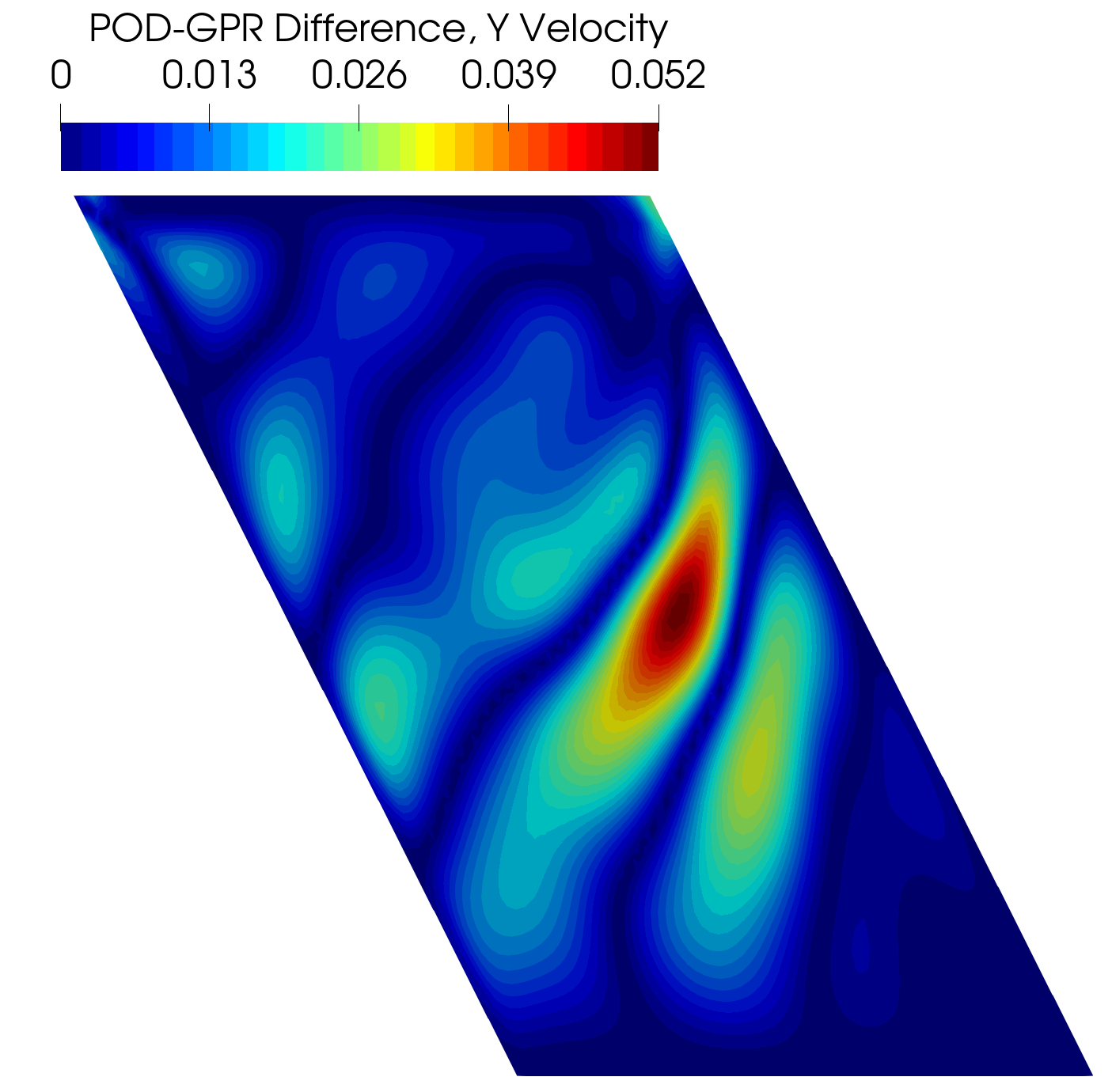}} 
\hspace{0.2cm}
\subfigure[$\bm{v}, \text{CAE-GPR Difference}$]{\includegraphics[scale=0.105]{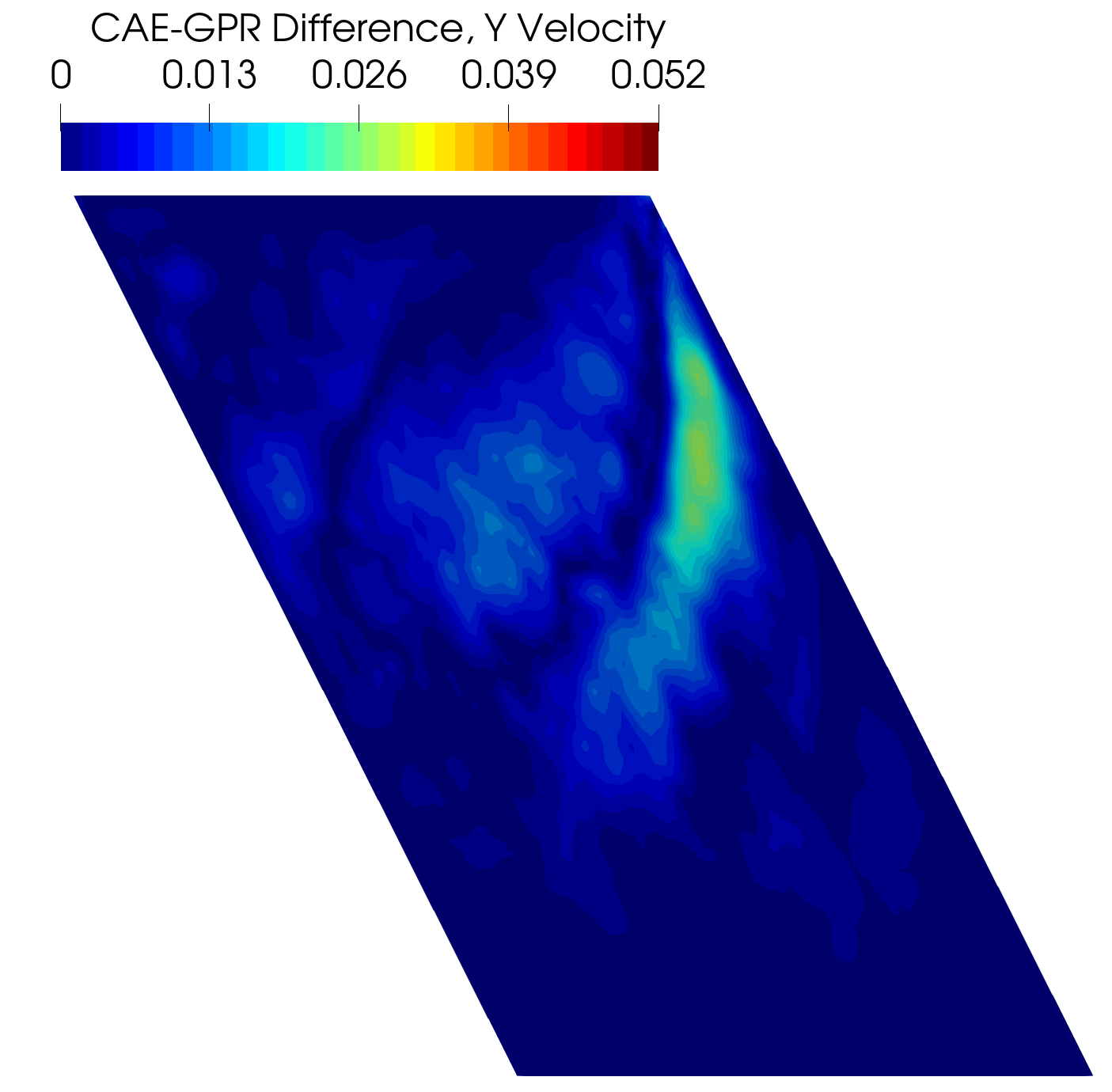}}
\caption{ROM comparison of $u$ and $v$ at $\bm{\mu} = (1.167, 1.997, -0.4665, 555.5)$, with $k = 5$ for CAE-GPR and $k = 35$ for POD-GPR.}
\label{fig:contour_322}
\end{figure}

\begin{figure}[!htbp]
\centering
\subfigure{\includegraphics[width=0.49\textwidth]{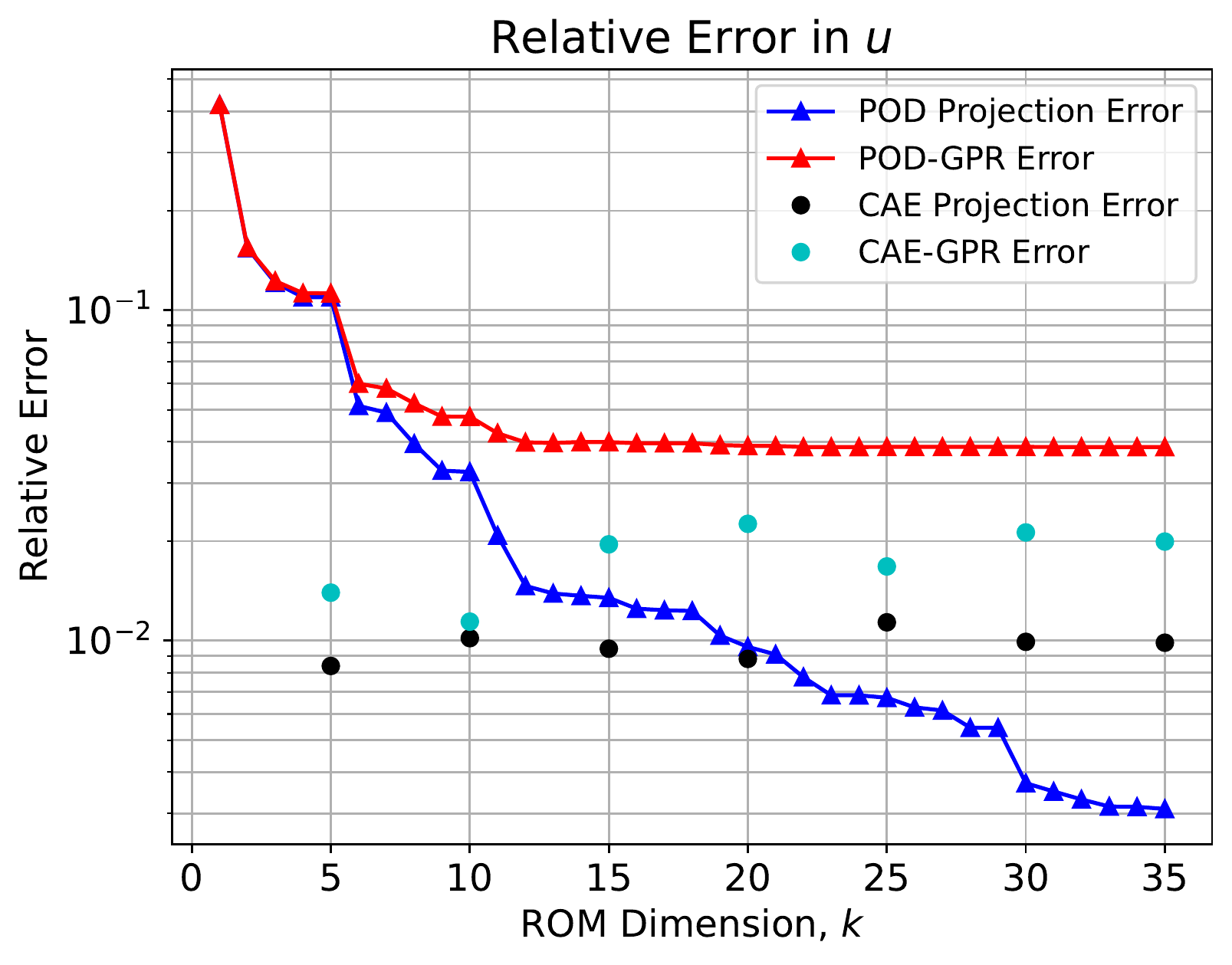}}
\subfigure{\includegraphics[width=0.49\textwidth]{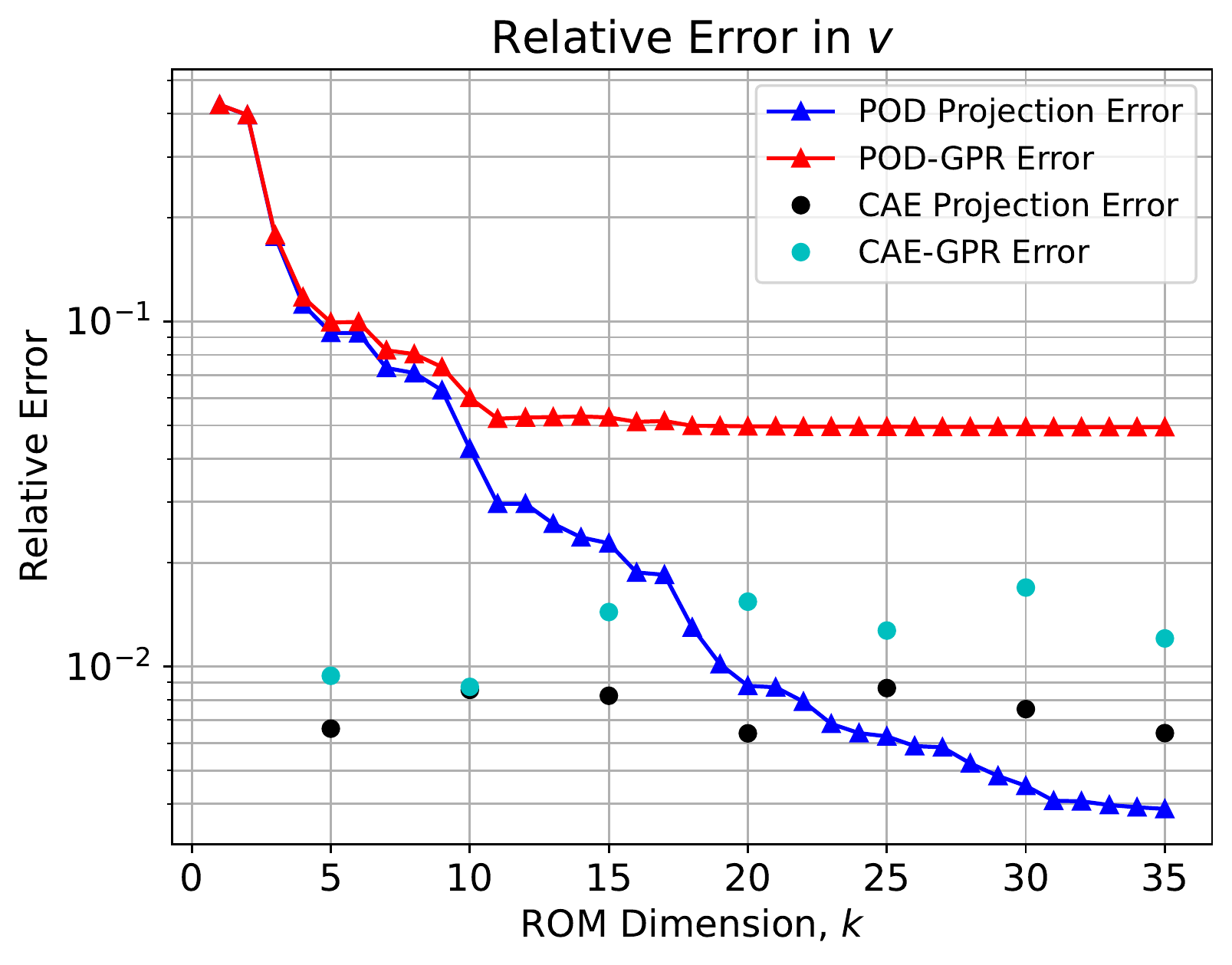}}
\caption{Plots of the prediction and projection errors in $u$ and $v$ for both ROMs at $\bm{\mu} = (1.963, 1.789, 0.5890, 308.5)$ at different values of $k$.} 
\label{fig:err_346}
\end{figure}

\begin{figure}[!t]
\centering
\subfigure[$\bm{u}, \text{Ground Truth}$]{\includegraphics[scale=0.112]{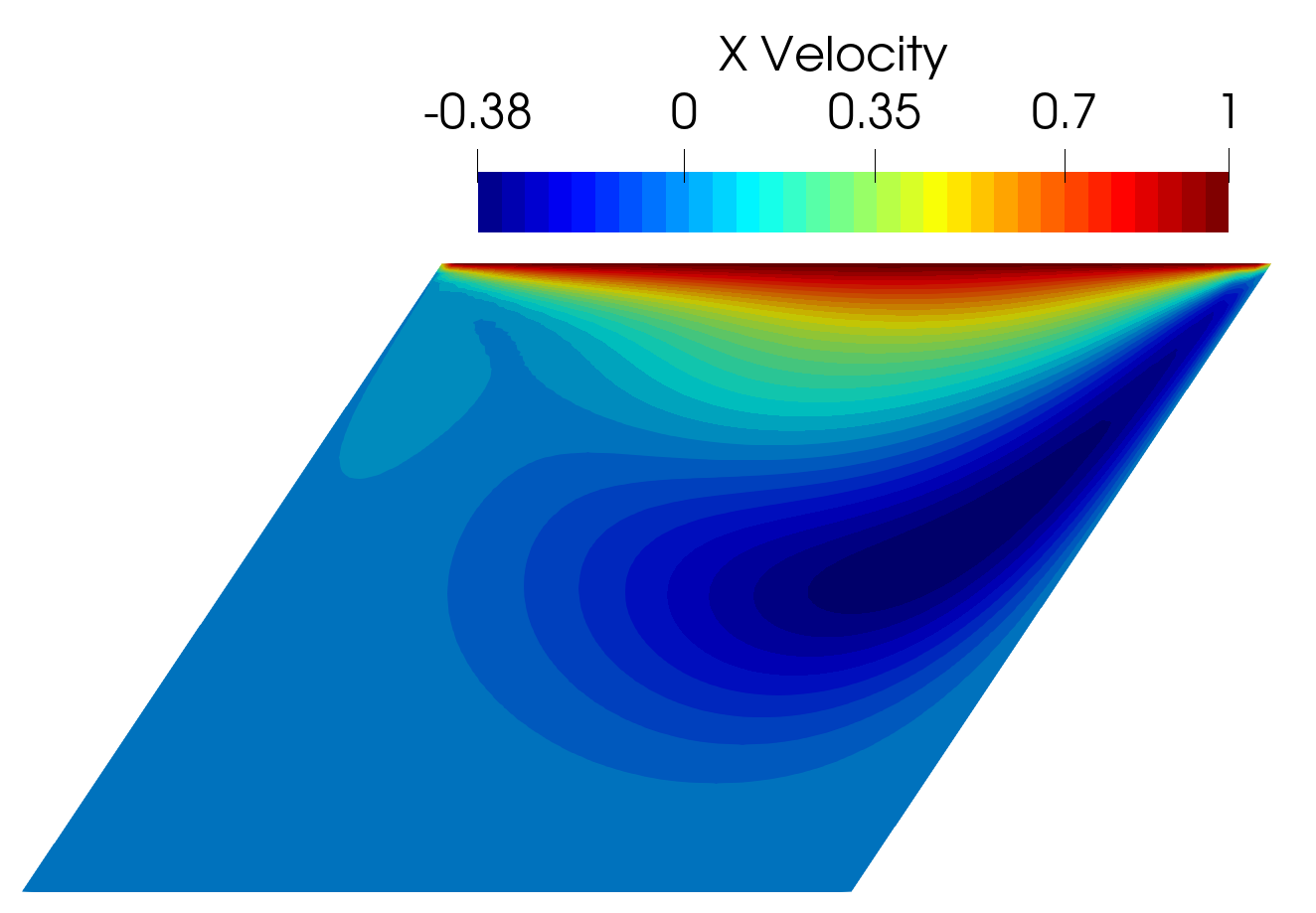}}
\hspace{0.2cm}
\subfigure[$\bm{u}, \text{POD-GPR Difference}$]{\includegraphics[scale=0.112]{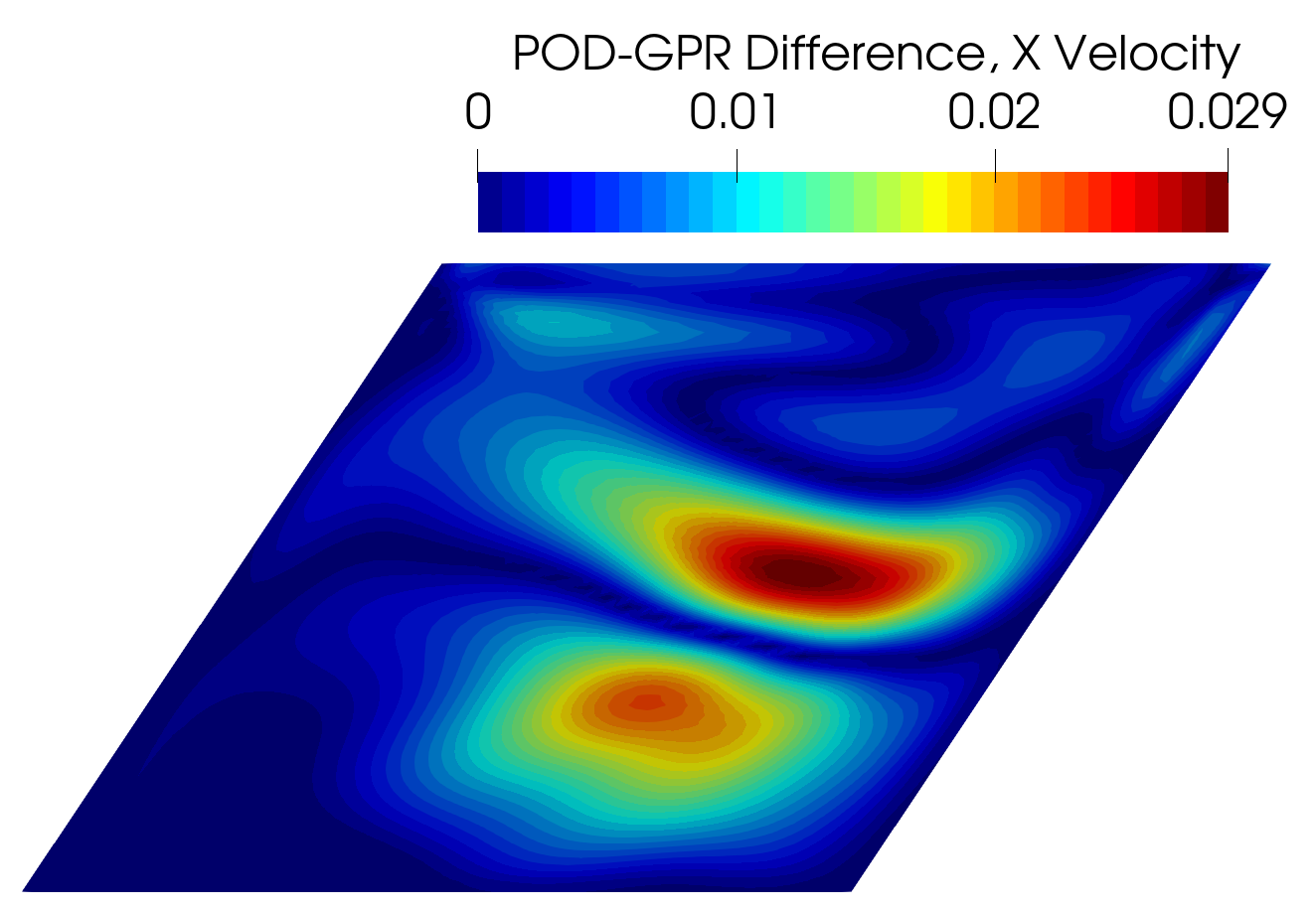}} 
\hspace{0.2cm}
\subfigure[$\bm{u}, \text{CAE-GPR Difference}$]{\includegraphics[scale=0.112]{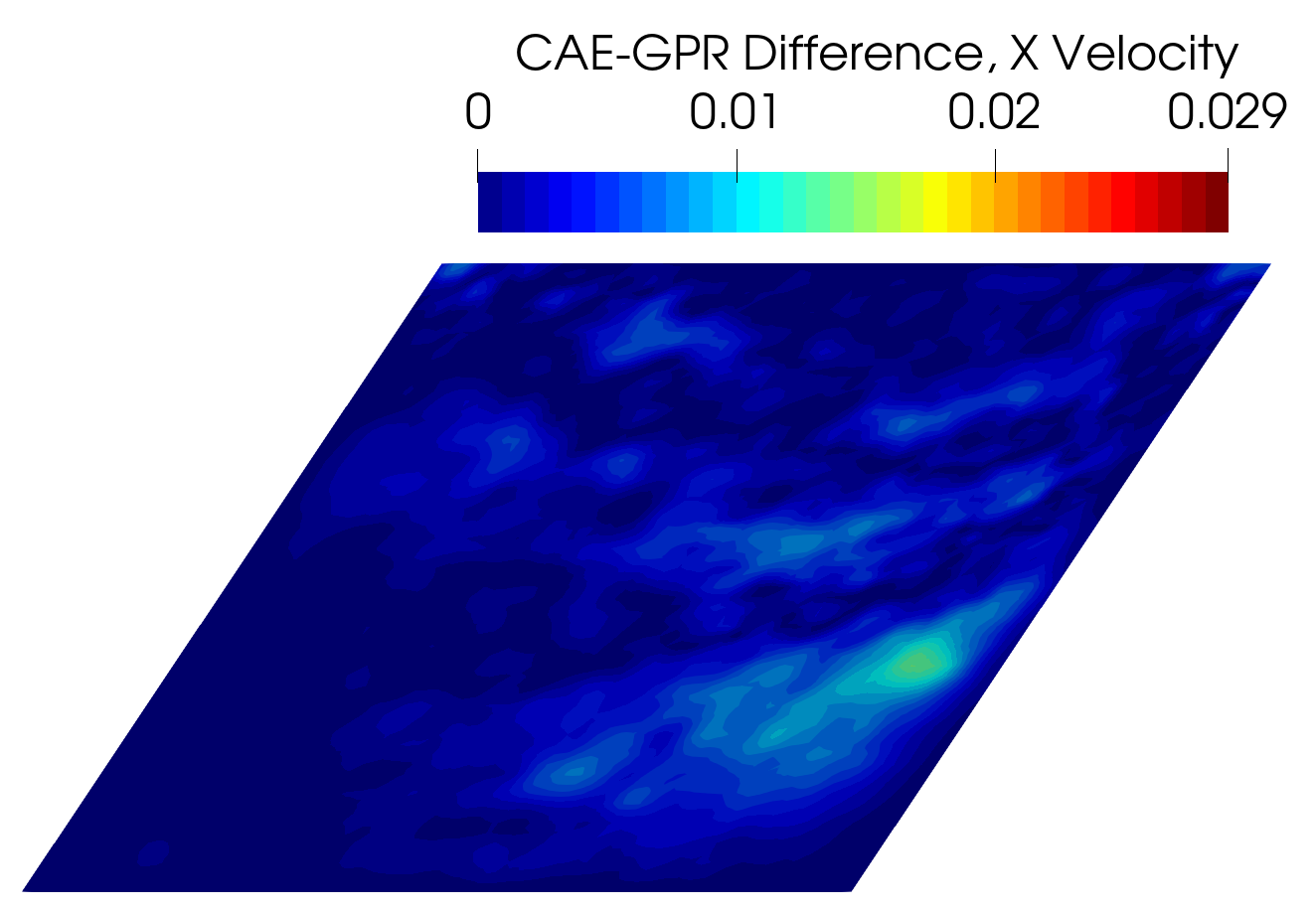}}
\subfigure[$\bm{v}, \text{Ground Truth}$]{\includegraphics[scale=0.112]{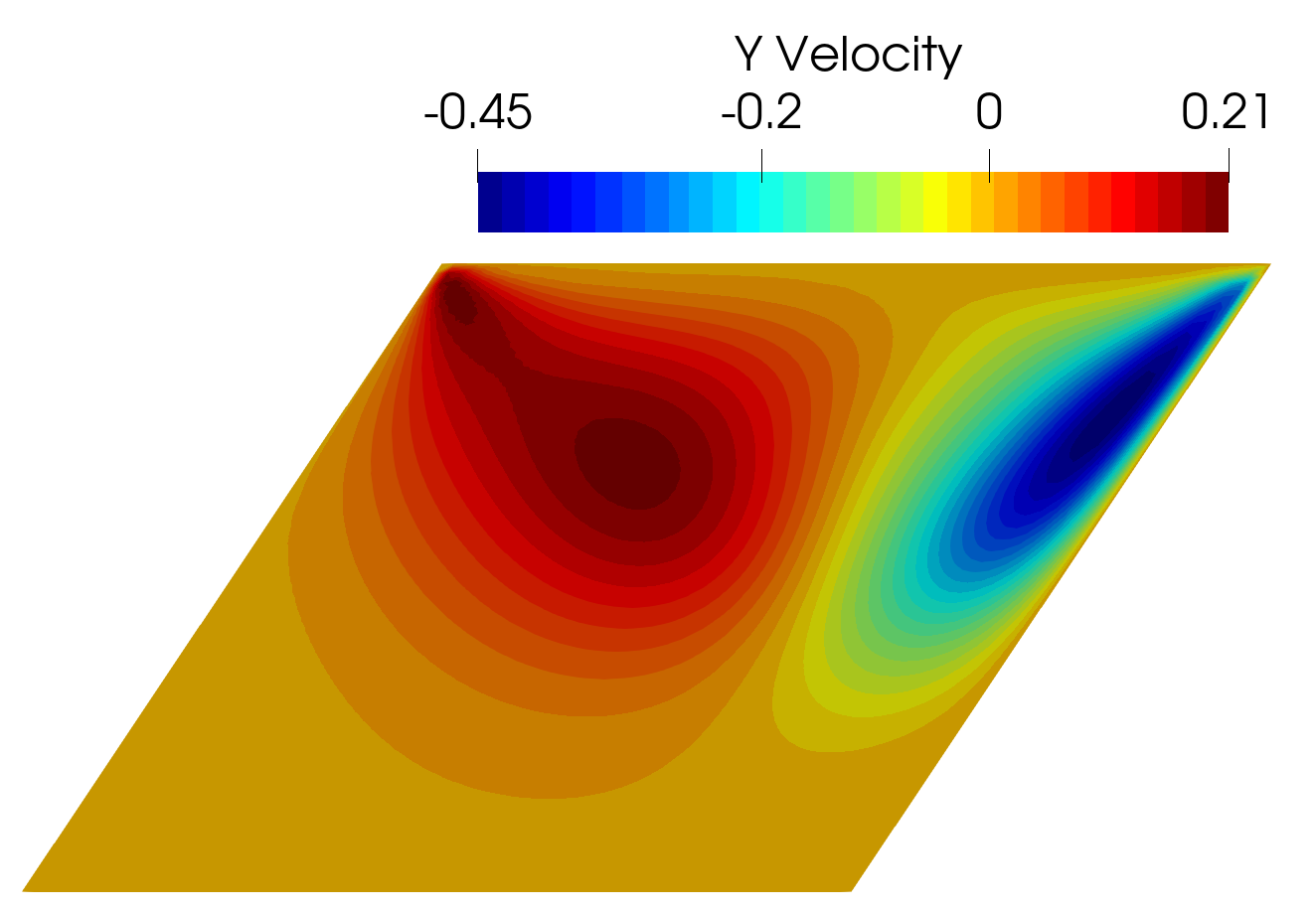}}
\hspace{0.2cm}
\subfigure[$\bm{v}, \text{POD-GPR Difference}$]{\includegraphics[scale=0.112]{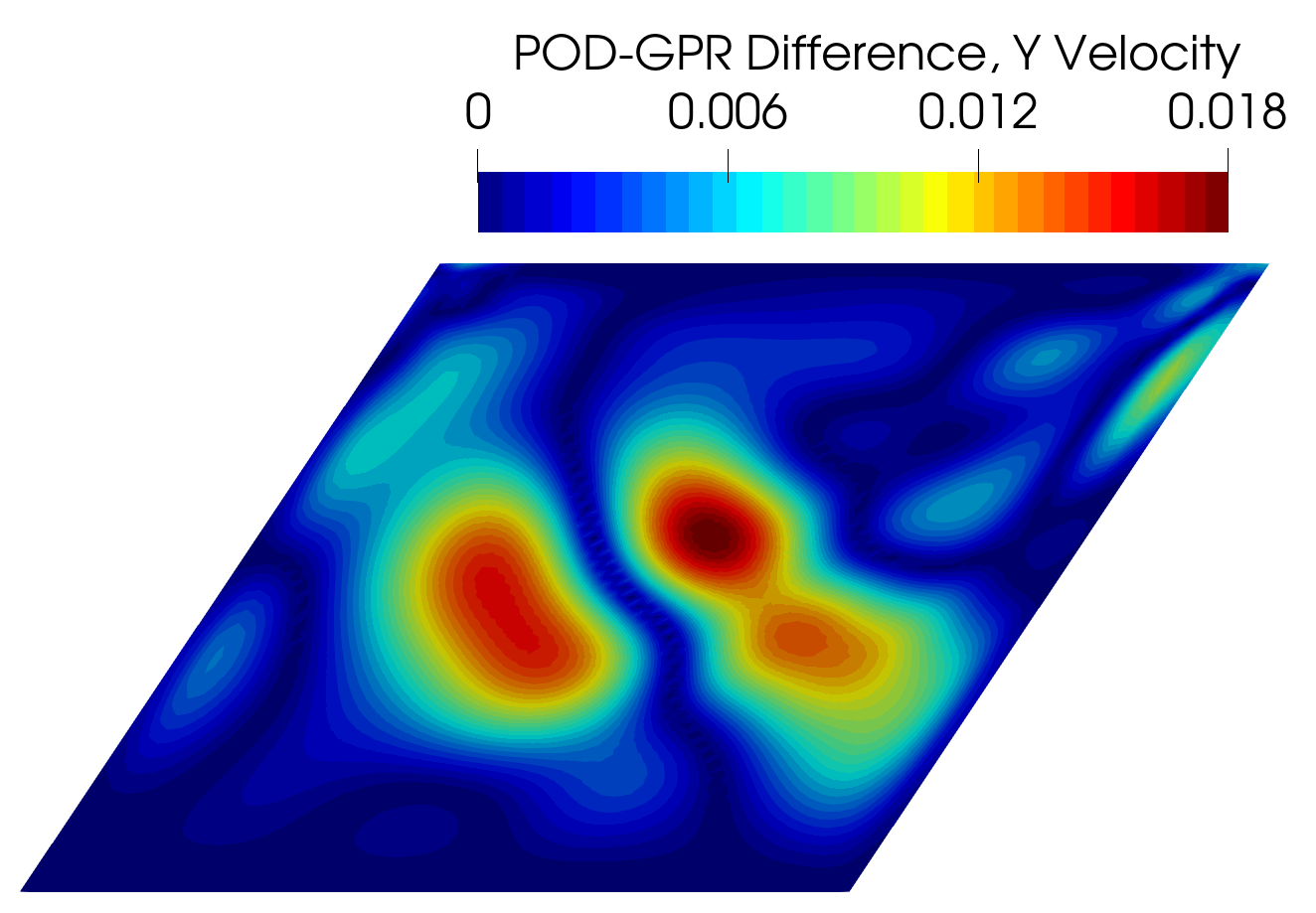}} 
\hspace{0.2cm}
\subfigure[$\bm{v}, \text{CAE-GPR Difference}$]{\includegraphics[scale=0.112]{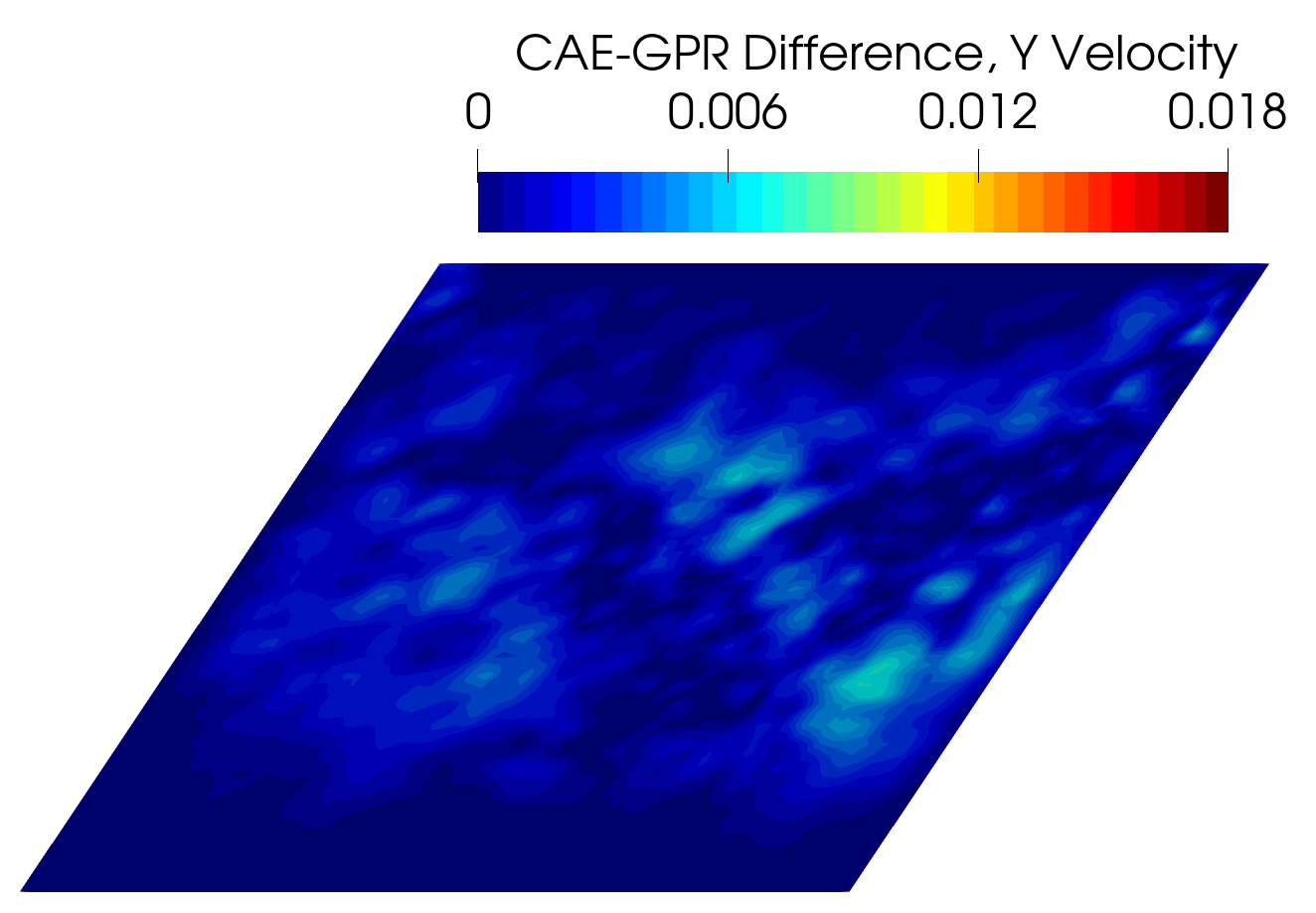}}
\caption{ROM comparison of $u$ and $v$ at $\bm{\mu} = (1.963, 1.789, 0.5890, 308.5)$, with $k = 5$ for CAE-GPR and $k = 35$ for POD-GPR.}
\label{fig:contour_346}
\end{figure}

Figure~\ref{fig:err_322} shows the relative error plot for both ROMs at a single design parameter $\bm{\mu}_{1} = (1.167, 1.997, -0.4665, 555.5)$, while Figure~\ref{fig:contour_322} shows the contour plots of the FOM as well as the absolute error plots for CAE-GPR at $k = 5$ and POD-GPR at $k = 35$. Results at another design parameter instance at $\bm{\mu}_{2} = (1.963, 1.789, 0.5890, 308.5)$ are shown in Figures~\ref{fig:err_346} and~\ref{fig:contour_346}. The generalized results from the cross-validation also hold here; for a greater projection error, CAE-GPR provides a lower prediction error. It is also shown at these design parameters that the prediction error curve of POD-GPR flattens out. There is more volatility in both the projection and prediction errors for CAE-GPR, although POD-GPR still never outperforms it in predicting $u$ and $v$. For $k = 10$ at both design parameters, the difference in the projection and prediction errors is very small, and almost 0 for $u$ at $\bm{\mu}_{1}$ and $v$ at $\bm{\mu}_{2}$. This is similar to the behavior exhibited by POD-GPR for low values of $k$, showing that CAE-GPR is also capable of producing highly accurate estimates of the expansion coefficients. The error contours at these design parameters highlight the increased predictive performance given by CAE-GPR over POD-GPR. While the error contours given by POD-GPR exhibit distinct bands of high error, the contours produced by CAE-GPR are generally more uniform and dispersed throughout the domain. At the chosen design parameters, there is a significant decrease in relative error; at $\bm{\mu}_{1}$, the percent decreases in relative error of $u$ and $v$ from POD-GPR to CAE-GPR are 42.4$\%$ and 49.8$\%$ respectively, while at $\bm{\mu}_{2}$ they are 48.3$\%$ and 75.6$\%$ respectively.

\section{Conclusion}
This work presents a non-intrusive reduced-order model framework utilizing nonlinear trial manifolds through the use of convolutional autoencoders. A deep learning approach, CAEs learn efficient low-dimensional representations of data through an encoder and decoder connected by a code.  CAEs provide a set of expansion coefficients through the low-dimensional code, similar to the basis coefficients provided by POD-based methods. A nonlinear relationship exists between the expansion coefficients and full-order states when using autoencoders, in contrast with the linear relationship when using POD. Given steady-state solutions of PDEs parameterized by a set of design parameters, Gaussian process regression can be used to approximate the expansion coefficients at unseen points in the design space for both approaches (referred to as CAE-GPR and POD-GPR). CAE-GPR involves a more expensive offline stage due to the high computational cost associated with training deep neural networks and requires that solutions have appropriate spatial arrangement. 

When applied to a geometrically and physically parameterized lid-driven cavity problem solved using the steady incompressible Navier-Stokes equations, it is shown that CAE-GPR offers higher performance in predicting the components of the velocity field when compared to POD-GPR over a range of ROM dimension. The low-dimensional code provided by CAEs is shown to be more easily interpolated than the basis coefficients obtained from POD. For a greater manifold projection error, CAE-GPR provides lower error in predictions of full-order states. It is also shown that CAE-GPR has the ability to provide highly accurate estimates of the expansion coefficients, providing prediction errors that are very close to projection errors. Although previous works~\cite{mrosek2021variational, kadeethum2022non} have shown that the autoencoders do not offer a remarkable advantage in performance over POD for some problems, highly non-linear problems such as the lid-driven cavity problem presented benefit significantly from the use of deep learning for ROM construction. Future work will extend this ROM framework to larger problems where the spatial arrangement of full-order states is not uniform and investigate constructing nonlinear trial manifolds using variational autoencoders (VAEs), which have been shown to provide a more interpretable low-dimensional code.

\begin{appendices}
\section{Convolutional autoencoder architecture}
The CAE architecture used for the lid-driven cavity problem ROM is listed in Table~\ref{tab:cae}. Zero padding is used for all convolutional and max-pooling layers. All leaky ReLU activation functions use a value of $\alpha = 0.25$. There are a relatively small number of convolutional and pooling layers in the network; we found that adding more of them did not improve the network performance, although their absence (using an MLP) causes a large decrease in performance. Compared to images which can be very noisy, the physical states in the presented problem vary smoothly, and fewer convolutional layers are required to learn features. We also found that having a fully connected layer on either side of the code was important for network performance, although this does drastically increase the number of network parameters. This makes the use of early stopping as a regularization method important to ensure that the network does not overfit. The learning task at hand requires that reconstructions of data be highly accurate, and networks with more parameters allow for more robust functional relationships to arise.
% \usepackage{booktabs}
\begin{table}[!htbp]
\centering
\begin{tabular}{@{}llllll@{}}
\toprule
Layer                   & Number of Filters & Kernel Size & Stride Value & Activation Function & Size of Output \\ \midrule
Input                   &                   &             &              &                     & 64 $\times$ 64 $\times$ 2    \\
Convolutional           & 64                & 3 $\times$ 3       & 1 $\times$ 1        & Leaky ReLU          & 64 $\times$ 64 $\times$ 64   \\
Max-Pooling             &                   & 2 $\times$ 2       & 2 $\times$ 2        &                     & 32 $\times$ 32 $\times$ 64   \\
Convolutional           & 32                & 3 $\times$ 3       & 1 $\times$ 1        & Leaky ReLU          & 32 $\times$ 32 $\times$ 32   \\
Max-Pooling             &                   & 2 $\times$ 2       & 2 $\times$ 2        &                     & 16 $\times$ 16 $\times$ 32   \\
Reshape                 &                   &             &              &                     & 8192           \\
Fully Connected         &                   &             &              & Leaky ReLU          & 128            \\
Code                    &                   &             &              & Leaky ReLU          & $k$              \\
Fully Connected         &                   &             &              & Leaky ReLU          & 128            \\
Fully Connected         &                   &             &              & Leaky ReLU          & 8192           \\
Reshape                 &                   &             &              &                     & 16 $\times$ 16 $\times$ 32   \\
Convolutional Transpose & 32                & 3 $\times$ 3       & 2 $\times$ 2        & Leaky ReLU          & 32 $\times$ 32 $\times$ 32   \\
Convolutional Transpose & 64                & 3 $\times$ 3       & 2 $\times$ 2        & Leaky ReLU          & 64 $\times$ 64 $\times$ 64   \\
Convolutional Transpose & 2                 & 3 $\times$ 3       & 1 $\times$ 1        & Sigmoid             & 64 $\times$ 64 $\times$ 2  \\ \bottomrule 
\end{tabular}
\caption{Detailed convolutional autoencoder (CAE) architecture used for the lid-driven cavity problem ROM.}
\label{tab:cae}
\end{table}

\section{Convolutional autoencoder training}

\begin{figure}[!htbp]
\centering
\subfigure{\includegraphics[width=0.47\textwidth]{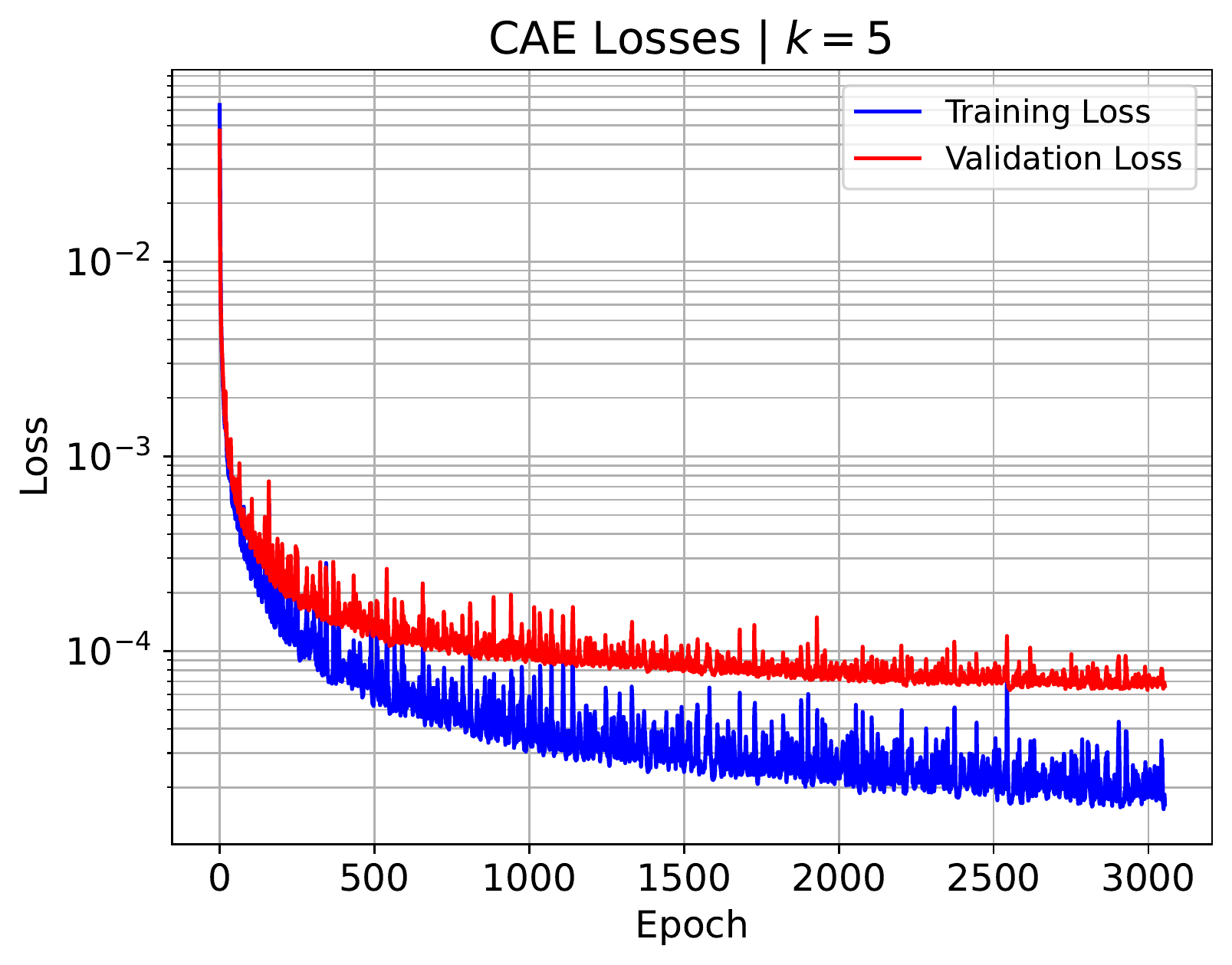}}
\subfigure{\includegraphics[width=0.47\textwidth]{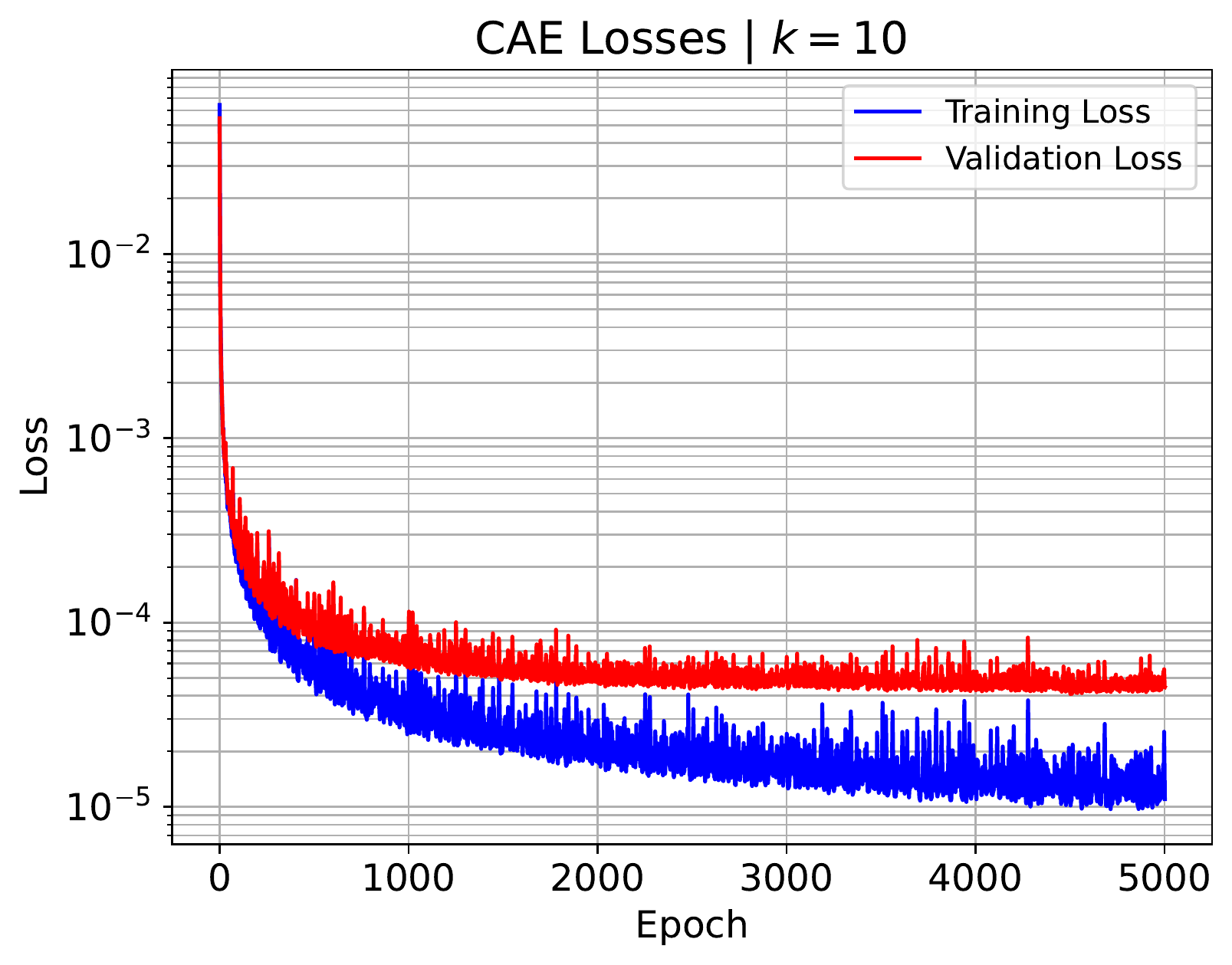}}
\subfigure{\includegraphics[width=0.47\textwidth]{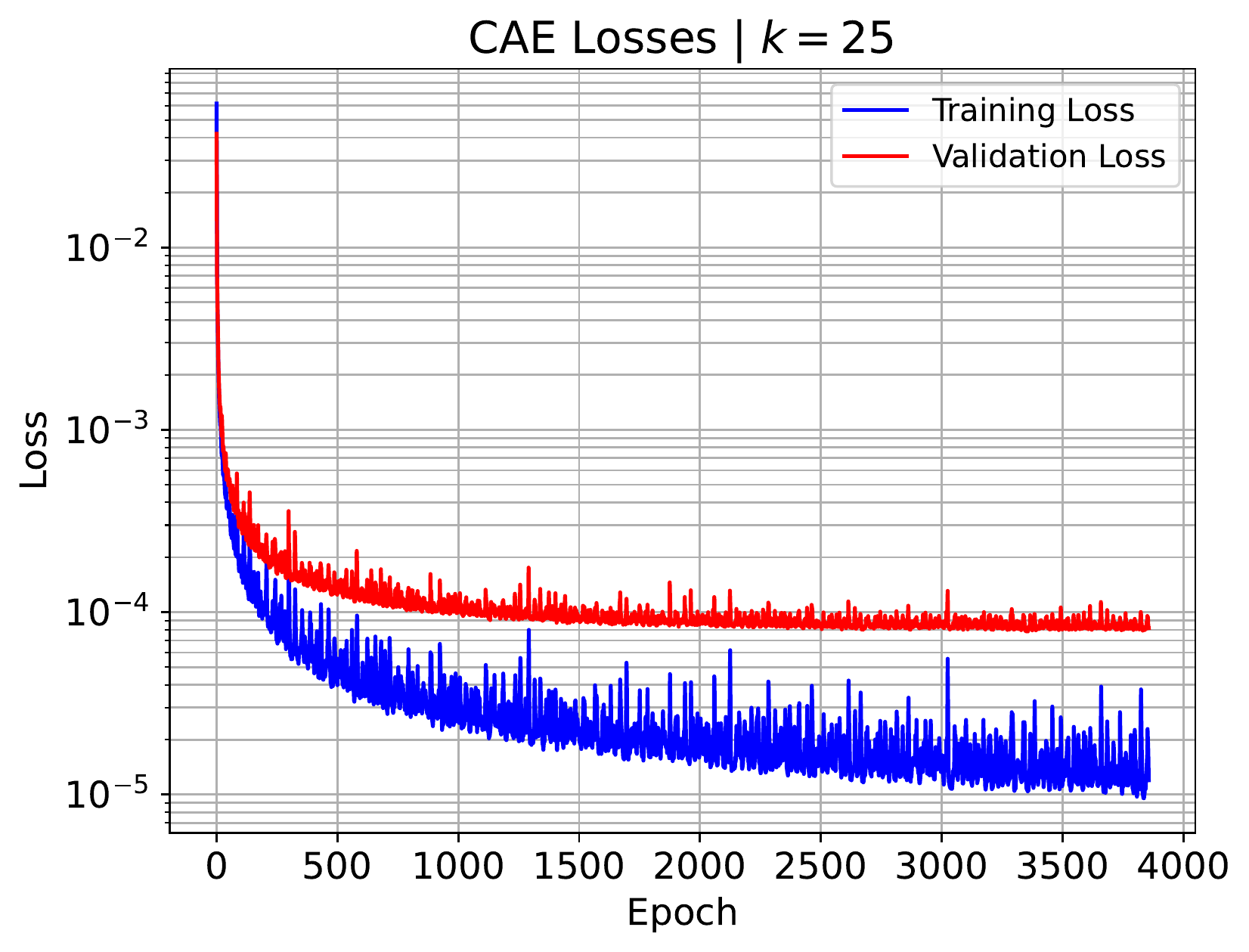}}
\subfigure{\includegraphics[width=0.47\textwidth]{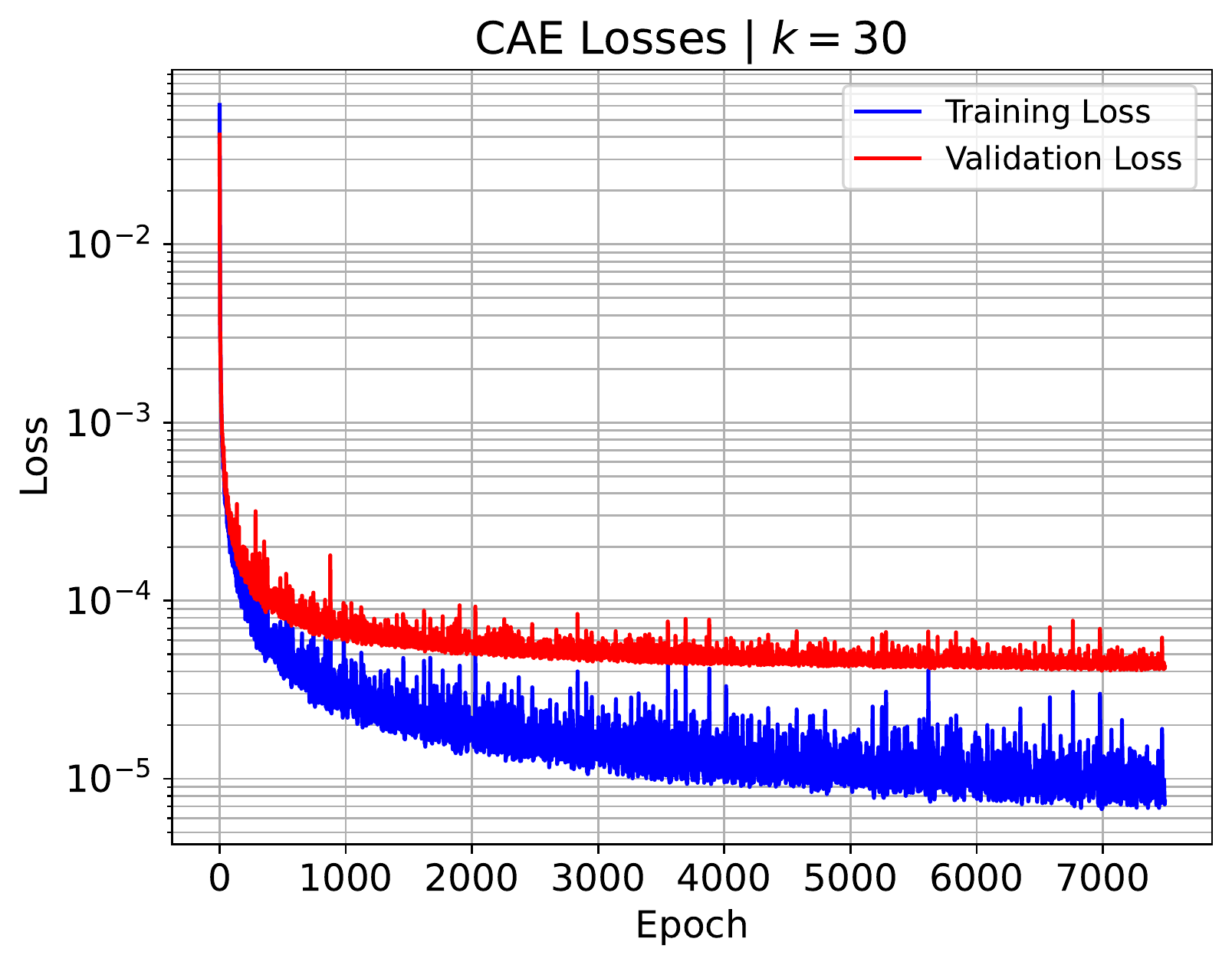}}
\caption{Plots of the training and validation losses at different ROM dimensions $k$ for a selected fold of the training and validation data.} 
\label{fig:cae_train}
\end{figure}

Figure~\ref{fig:cae_train} shows the training and validation losses against the number of epochs for selected folds of the training and validation data for $k = 5, 10, 25, 30.$ Early stopping is used as a regularization method, and the validation loss fails to drop for 500 epochs well before the maximum number of 7500 epochs at $k = 5, 10, 25$. At $k = 30$, training stops after 7491 epochs. While the training loss continues to decline slowly in all of the plots, the validation loss shows asymptotic behavior. By monitoring the validation loss and using early stopping, the network is prevented from overfitting the training data. Training is performed on an NVIDIA TITAN RTX GPU. The average wall time and number of epochs for training the CAE over all of the data folds is shown in Table~\ref{tab:costs}; in general, increasing $k$ leads to higher computational costs. Both the number of trainable parameters and capacity of the network to learn are affected by the size of the code.

% \usepackage{booktabs}
\begin{table}[!htbp]
\centering
\begin{tabular}{@{}lll@{}}
\toprule
ROM Dimension, $k$ & Average Wall Time (s) & Average Number of Epochs \\ \midrule
5             & 884                   & 3533                     \\
10            & 883                   & 3564                     \\
15            & 867                   & 3441                     \\
20            & 1123                  & 4465                     \\
25            & 971                   & 3865                     \\
30            & 1283                  & 4675                     \\
35            & 1301                  & 4519                     \\ \bottomrule 

\end{tabular}
\caption{Average computational costs over all data folds for training the CAE.}
\label{tab:costs}
\end{table}

\end{appendices}
\clearpage
\bibliographystyle{unsrt}  
\bibliography{references}

\end{document}